\documentclass[11pt,reqno]{amsart}
\usepackage{latexsym,amsmath,amssymb,amsfonts,amsthm,verbatim}
\usepackage[utf8]{inputenc}
\usepackage{color}

\let \ssection=\section
\renewcommand{\section}{\setcounter{equation}{0}\ssection}
\newtheorem{Th}{Theorem}[section]
\newtheorem{theorem}[Th]{Theorem}
\newtheorem{corollary}[Th]{Corollary}
\newtheorem{proposition}[Th]{Proposition}
\newtheorem{lemma}[Th]{Lemma}
\theoremstyle{definition}
\newtheorem{definition}[Th]{Definition}

\theoremstyle{remark}
\newtheorem{remark}[Th]{Remark}
\newcommand{\R}{\mathbb R}
\newcommand{\E}{\mathbb E}
\newcommand{\Ek}{\mathbb E_k}
\newcommand{\F}{\mathcal F}
\newcommand{\N}{\mathcal N}
\newcommand{\rv}{\right|}
\newcommand{\lv}{\left|}
\newcommand{\lp}{\left (}
\newcommand{\rp}{\right )}
\newcommand{\lc}{\left [}
\newcommand{\rc}{\right ]}
\newcommand{\e}{\epsilon}

\newcommand{\dBr}{d\overleftarrow{B}_r}
\newcommand{\dBs}{d\overleftarrow{B}_s}
\newcommand{\imImM}{^{i-1,I,m,M}}
\newcommand{\YNt}[1]{Y^N_{ t_{#1}}}
\newcommand{\XNt}[1]{X^N_{t_{#1}}}
\newcommand{\ZNt}[1]{Z^N_{ t_{#1}}}

\newcommand{\unh}{\frac{1}{h}}
\newcommand{\DeltaW}[1]{\Delta W_{#1}}
\newcommand{\DeltaB}[1]{\overleftarrow{\Delta} B_{#1}}
\newcommand{\Ftk}[1][k]{\mathcal F_{t_{#1}}}
\newcommand{\XYZr}{\lp X_r,Y_r,Z_r\rp}

\newcommand{\XYZs}{\lp X_s,Y_s,Z_s\rp}
\newcommand{\XYs}{\lp X_s,Y_s\rp}
\newcommand{\XYZNt}[1]{\lp \XNt{#1},\YNt{#1},\ZNt{#1}\rp}
\newcommand{\undeux}{\frac{1}{2}}
\newcommand{\YNIt}[2]{Y^{N,#1,I}_{t_{#2}}}

\newcommand{\ay}[2]{\alpha^{#1,I}_{#2}}
\newcommand{\az}[2]{\beta^{#1}_{#2}}
\newcommand{\py}[1]{p_{#1}}
\newcommand{\pz}[1]{p_{#1}}
\newcommand{\tk}{_{t_k}}
\newcommand{\tkk}{_{t_{k+1}}}
\renewcommand{\P}{\mathcal P}
\newcommand{\NooI}{^{N,\infty,I}}
\newcommand{\NiI}{^{N,i,I}}
\newcommand{\NII}{^{N,I,I}}
\newcommand{\pY}{Y^{\pi}}
\newcommand{\pZ}{Z^{\pi}}
\newcommand{\uZ}{Z^{\pi,1}}
\newcommand{\tkm}{_{t_{k-1}}}

\newcommand{\VMk}{V^M_k}
\newcommand{\AMk}{A^M_k}

\newcommand{\aIIMkk}{\alpha^{I,I,M}_{k+1}}

\newcommand{\aiIMk}{\alpha^{\infty,I,M}_{k}}

\newcommand{\biIMk}{\beta^{\infty,I,M}_{k}}
\newcommand{\qmk}{p^m_{k}}
\newcommand{\pmk}{p^m_{k}}
\newcommand{\pmkk}{p^m_{k+1}}
\newcommand{\XNm}{X^{N,m}}
\newcommand{\DeltaBmk}{\overleftarrow{\Delta} B^m_k}

\newcommand{\NIIM}{^{N,I,I,M}}
\newcommand{\DeltaWm}[1]{\Delta W^m_{#1}}
\newcommand{\DeltaBm}[1]{\overleftarrow\Delta B^m_{#1}}
\newcommand{\IIM}{^{I,I,M}}
\newcommand{\Mk}{^M_k}
\newcommand{\unM}{\frac{1}{M}}
\newcommand{\summ}{\sum_{m=1}^M}
\newcommand{\mk}{^m_k}
\newcommand{\gXYr}{g \lp X_r,Y_r\rp}
\newcommand{\gXYpk}{g\lp X^N\tk, Y^{\pi}\tk\rp}
\newcommand{\NinftyI}{^{N,\infty,I}}
\newcommand{\NI}{^{N,I}}
\newcommand{\Nk}{^N_k}
\newcommand{\wz}{\widehat\zeta\Nk}
\renewcommand{\wr}{\widehat\rho\Nk}
\newcommand{\kk}{_{k+1}}
\newcommand{\Nj}{^N_j}
\newcommand{\AMj}{A^M_j}
\newcommand{\II}{^{I,I}}
\newcommand{\iIM}{^{i,I,M}}
\newcommand{\imIM}{^{i-1,I,M}}
\newcommand{\iiIM}{^{i+1,I,M}}
\newcommand{\Nm}{^{N,m}}
\newcommand{\inftyIM}{^{\infty,I,M}}
\newcommand{\inftyImM}{^{\infty,I,m,M}}
\newcommand{\wrNm}{\widehat \rho^{N,m}}
\newcommand{\mkk}{^{m}_{k+1}}
\newcommand{\inftyI}{^{\infty,I}}
\renewcommand{\k}{_k}
\newcommand{\VMkm}{\lp V^M_k\rp^{-1}}
\newcommand{\wrN}{\widehat \rho^N}
\newcommand{\iI}{^{i,I}}
\newcommand{\iImM}{^{i,I,m,M}}

\newcommand{\rvd}{\right |^2}
\newcommand{\VM}{V^M}
\newcommand{\EM}{\mathfrak E_M}
\newcommand{\B}{\mathfrak B}

\newcommand{\lb}{\left\lbrace}
\newcommand{\rb}{\right\rbrace}

\newcommand{\tinitial}{0}
\newcommand{\PP}{\mathbb P}
\newcommand{\Fchapeau}{\widehat{\mathcal F}}

\newcommand{\sumLambda}{}
\newcommand{\gl}{g}
\newcommand{\dB}[1]{d\overleftarrow{B}_{#1}}

\newcommand{\uik}{u_i\lp X^N_{t_k}\rp}
\newcommand{\uiNk}{u_{i_N}\lp X^N_{t_k}\rp}

\newcommand{\vik}[1]{v_i\lp#1\rp}
\newcommand{\DeltahB}[1]{\frac{\overleftarrow\Delta B_{#1}}{\sqrt h}}

\newcommand{\iN}{{i_N}}
\newcommand{\iNm}{{i_{N-1}}}
\newcommand{\ikk}{{i_{k+1}}}
\newcommand{\viNmk}[1]{v_{i_{N-1}}\lp#1\rp}
\newcommand{\vikkk}[1]{v_{i_{k+1}}\lp#1\rp}
\newcommand{\gXYNkk}{g\lp X^N\tkk,Y^N\tkk\rp}
\newcommand{\gXYNk}{g\lp X^N\tk,Y^N\tk\rp}
\newcommand{\Ekm}{\E_{k-1}}
\newcommand{\XYZpikm}{\lp X^N\tkm, \pY\tkm, \uZ\tkm\rp}

\newcommand{\km}{_{k-1}}
\newcommand{\tC}{\widetilde C}
\newcommand{\bC}{\overline{C}}
\newcommand{\gTNIkk}{g\lp X^N_{t_{k+1}},Y^{N,I,I}\tkk\rp}
\newcommand{\Yinf}[1]{Y^{N,\infty,I}_{t_{#1}}}
\newcommand{\XYZinf}[1]{\lp X^N_{t_{#1}},\Yinf{#1},Z^{N,I}_{t_{#1}}\rp}
\newcommand{\IM}{^{I,M}}

\begin{document}
\title[Regression Monte-Carlo for BDSDE]
{A regression Monte-Carlo method for
Backward Doubly Stochastic Differential
Equations}

\author{Omar Aboura}
\date{\today}

\address{ 
  SAMM (EA 4543),
 Universit\'e Paris 1 Panth\'eon Sorbonne,
90 Rue de Tolbiac, 75634 Paris Cedex France }
 \email{omar.aboura@malix.univ-paris1.fr}
\maketitle
\begin{abstract}
This paper extends the idea of E.Gobet, J.P.Lemor and X.Warin 
from the setting of Backward Stochastic Differential Equations to that
of Backward Doubly Stochastic Differential equations. We propose
some numerical approximation scheme of these equations introduced by
E.Pardoux and S.Peng.
\end{abstract}

\section{Introduction} 

Since the pioneering work of E.~Pardoux and S.~Peng  \cite{PPquasilinear}, 
 backward stochastic differential equations (BSDEs) have been intensively studied during the two last decades. 
Indeed, this notion has been a very useful tool to study problems in many areas, such as
mathematical finance, stochastic control, partial differential equations; see e.g.  \cite{MaYong}
where many applications are described.
Discretization schemes for BSDEs have been studied by several authors. 
The first papers on this  topic are that of V.Bally \cite{BallyDiscretization}  and D.Chevance \cite{ChevanceThesis}.
In his thesis, J.Zhang made an interesting contribution which was the starting point of intense study among,
which the works of   B.~Bouchard and N.Touzi \cite{BouchardTouzi}, E.Gobet, J.P. Lemor and
X. Warin\cite{Gobet},...
The notion of BSDE has been generalized by E.~Pardoux and S.~Peng  \cite{PPspde} to that of 
 Backward Doubly Stochastic Differential Equation (BDSDE) as follows. Let  $(\Omega,\F,\PP)$ be
a  probability space,  $T$ denote some fixed terminal time which will be used throughout the paper, 
$\lp W_t\rp_{\tinitial\leq t\leq T}$ and $\lp B_t\rp_{\tinitial\leq t\leq T}$ be two independent standard
 Brownian motions defined on  $(\Omega,\F,\PP)$ and  with values in $\R$. 
 On this space  we will deal with the following families of $\sigma$-algebras:
\begin{equation}\label{defF}
 \F_t:=\F^W_{\tinitial,t}\vee\F_{t,T}^B\vee\N,\quad
\Fchapeau_t:=\F^W_{\tinitial,t}\vee\F_{\tinitial,T}^B\vee\N, \quad 
{\mathcal H}_t=\F^W_{\tinitial,T}\vee\F_{t,T}^B\vee\N,
\end{equation}
where $\F_{t,T}^B:=\sigma\lp B_r-B_t; t\leq r\leq T\rp$, 
$\F^W_{\tinitial,t}:=\sigma\lp W_r;\tinitial\leq r\leq t \rp$ and
$\N$ denotes the class of $\PP$ null sets. We remark that  $(\Fchapeau_t)$ is a filtration, 
$({\mathcal H}_t)$ is a decreasing family of $\sigma$-albegras,  while 
$(\F_t)$ is neither  increasing nor decreasing. 
Given an initial condition  $x\in \R$, let $(X_t)$ be the diffusion process
 defined by
\begin{equation}\label{defX}
X_t = x+\int_{0}^tb\lp X_s\rp ds+\int_{0}^t\sigma\lp X_s\rp dW_s .
\end{equation}
Let  $\xi\in L^2(\Omega)$ be an $\R$-valued, $\F_T$-measurable random variable, 
$f$ and $g$ be  regular enough coefficients; 
consider the BDSDE 
defined as follows: 
\begin{align} 
 Y_t\; = \; \xi & 
+\int_t^Tf\lp s,X_s,Y_s, Z_s\rp ds
+\int_t^T\gl\lp s,X_s,Y_s,Z_s\rp \dBs -\int_t^TZ_sdW_s.\label{BDSDE}
\end{align}
In this equation,  $dW$ is the forward stochastic
 integral and $\dB{}$ is the backward stochastic integral
 (we send the reader to \cite{NualartPardoux anticipating} for more  details on backward integration).
A solution to \eqref{BDSDE} is a  pair of  real-valued process $(Y_t,Z_t)$,
 such that $Y_t$ and $Z_t$ are $\F_t$-measurable for every $t\in [0,T]$, such that \eqref{BDSDE} is satisfied  and
\begin{equation} \label{bndYZ}
 \E\Big(\sup_{0\leq s\leq T} |Y_s|^2\Big) +  \E \int_0^T |Z_s|^2 ds<+\infty.
\end{equation}
In \cite{PPspde} Pardoux and Peng have proved that under some 
 Lipschitz property on $f$ and $g$ which will be stated later,
\eqref{BDSDE} has a unique solution  $(Y,Z)$. 
They also proved that 
\begin{align*}
Y_t= u\lp t,X_t, \lp\DeltaB{s}\rp_{t\leq s\leq T}\rp,\quad
Z_t=v\lp t,X_t, \lp\DeltaB{s}\rp_{t\leq s\leq T}\rp,
\end{align*}
for some Borel functions $u$ and $v$.

The time discretization of BDSDEs has been addressed in \cite{Aman} when
 the coefficient $g$ does not
depend on $Z$; see also \cite{Aboura} in the more general setting for $g$ which may also depend on $Z$ as
in \cite{PPspde}. Both papers follow Zhang's approach and provide a theoretical approximation only using a constant time
mesh.

In order to obtain a more tractable discretization which could be implemented, a natural idea is to see whether the methods
introduced in \cite{Gobet} can be extended from the framework of BSDEs to that more involved of BDSDEs ; 
this is the aim of this paper.

We use three consecutive steps, 
and each time we give a precise estimate of the corresponding 
error.
Thus, we start with a time discretization $(Y^N_{t_k}, Z^N_{t_k})$ 
with a constant time mesh $T/N$.
We can prove  that 
\begin{align*}
Y^N\tk=u^N\lp t_k,X^N\tk,\DeltaB{N-1},\dots,\DeltaB{k}\rp,\quad
Z^N\tk=v^N\lp t_k,X^N\tk,\DeltaB{N-1},\dots,\DeltaB{k}\rp,
\end{align*}
where for $k=1,\dots,N-1$, $t\k=kT/N$ and $\DeltaB{k}=B\tkk-B\tk$. 
Furthermore, if either $f=0$ or if the scheme is not implicit as in \cite{Aboura}
then we have the more precise description:
\begin{align*}
Y^N\tk= u^N_N\lp t_k,X^N\tk\rp
+\sum_{j=k}^{N-1}u^N_j\lp t_k,X^N\tk,\DeltaB{N-1},\dots,\DeltaB{j+1}\rp\DeltaB{j},\\
Z^N\tk= v^N_N\lp t_k,X^N\tk\rp
+\sum_{j=k}^{N-1}v^N_j\lp t_k,X^N\tk,\DeltaB{N-1},\dots,\DeltaB{j+1}\rp\DeltaB{j},
\end{align*}
with the convention that if $j+1>N-1$, $\lp \DeltaB{N-1},\dots,\DeltaB{j+1}\rp=\emptyset$.
The main time discretization result in this direction 
is Theorem \ref{theorem 1 Gobet}.
In order to have a numerical scheme, we use this decomposition and the ideas of E.Gobet, J.P.Lemor
and X.Warin \cite{Gobet}. Thus we introduce the following hypercubes, that is approximate random variables 
$u^N_j\lp t_k,X^N\tk,\DeltaB{N-1},\dots,\DeltaB{j+1}\rp\DeltaB{j}$ by their orthogonal projection on
some finite vector space generated by some bases $\lp u_j\rp$ and $\lp v_j\rp$ 
defined below. For $k=1,\dots,N$ we have\\
\begin{align*}
Y^N\tk
\approx&
\sum_{\iN}\E \lp Y^N\tk\uiNk \rp \uiNk\\
&+\sum_{j=k}^{N-1}
\sum_{\iN,\iNm,\dots,i_{j}}
\E\lp Y^N\tk\uiNk\viNmk{\DeltaB{N-1}}\dots\vikkk{\DeltaB{j+1}}\DeltahB{j}
\rp\\
&\qquad\qquad\qquad\qquad\uiNk\viNmk{\DeltaB{N-1}}\dots\vikkk{\DeltaB{j+1}}\DeltahB{j}.
\end{align*}

We use 
 a linear regression operator of the approximate solution. 
Thus, we at first use an orthogonal projection
 on a finite dimensional space ${\mathcal P}_k$. 
This space consists in linear
combinations of an orthonormal family of properly renormalized indicator functions of disjoint intervals composed either with
the diffusion $X$ or with increments of the Brownian motion $B$. 
 As in \cite{Gobet}, in order not to introduce error terms
worse that those due to the time discretization, 
we furtherore have to use a  Picard iteration scheme. 
The error due to this regression operator is estimated 
in Theorem \ref{th-2}.

Then the coefficients $(\alpha,\beta)$ of the decomposition
of the projection of $(Y^N_{t_k}, Z^N_{t_k})$ are shown to solve a regression minimization problem and are expressed in terms of expected values.
Note that a general regression approach has also been used 
by Bouchard and Touzi for BSDEs in \cite{BouchardTouzi}. 
Finally, the last step consists in replacing the minimization problem for the pair $(\alpha,\beta)$ in terms of expectations
by similar expressions described in terms of an average over a sample of size $M$ of the Brownian motions $W$ and $B$. Then, a proper
localization is needed to get an $L^2$ bound of the last error term. This requires another Picard iteration and the error term
due to this Monte Carlo method is described in Theorem \ref{theorem-stepIII}.  

A motivation to study BSDEs is that these equations are widely used 
in financial models,
so that having an efficient and fast numerical
methods is important.
As noted in \cite{PPspde}, BDSDEs are connected with stochastic partial differential equations and the discretization
of \eqref{equationYZ} is motivated by its link with the following SPDE:
\begin{eqnarray}
u(t,x)&=&\phi(x)+\int_t^T\Big(\mathcal L u(s,x)+f\lp s,x,u(s, x),\nabla u(s, x)\sigma (x)\rp \Big) ds\nonumber\\
&&
+\sumLambda\int_t^T\gl\lp s,x,u(s,x),\nabla u(s, x)\sigma (x)\rp \dB{s},\label{SPDE}
\end{eqnarray}
Discretizations of SPDEs are mainly based on PDE techniques, such as finite differences or finite elements methods. Another approach
for special equations is given by particle systems. We believe that this paper gives a third way to deal with this problem.
As usual, the presence of the gradient in the diffusion coefficient is the most difficult part to handle when dealing
with SPDEs.
Only few results are obtained in the classical discretization framework when PDE methods are extended to the stochastic case. 

Despite the fact that references \cite{Aman} and \cite{Aman2}
deal with a problem similar to that we address in section 3,
we have kept the results and proofs of this section.
Indeed, 
on one hand we study here an implicit scheme
as in \cite{Gobet} and wanted the paper to be self contained.
Furthermore,
because of measurability properties of $Y_0$ and $Y^{\pi}_0$, 
the statements and proofs of Theorem 3.6 in \cite{Aman} and Theorem 
4.6 in \cite{Aman2} are unclear and there is a gap in the corresponding
proofs because of similar measurability issues for $(Y_t)$ and $(Y^{\pi}_t)$.

The paper is organized as follows. Section 2 gives the main notations concerning the time discretization and the function basis.
Section 3 describes the time discretization and results similar to those in \cite{Aman} are proved in a more general
framwork. The fourth section describes the projection error. Finally section 5 studies the regression technique and
the corresponding  Monte Carlo method. Note that the presence of increments of the Brownian motion $B$, which drives the
backward stochastic integrals, requires some new arguments such as Lemma \ref{lem-avant-R} which is a key ingredient of the last error 
estimates. 
As usual $C$ denotes a constant which can change from line to line.
\section{Notations}

Let $\lp W_t,t\geq0\rp$ and $\lp B_t,t\geq0\rp$ 
be two mutually independent standard Brownian motions.
For each $x\in \R$, let $\lp X_t, Y_t,Z_t,t\in[0,T]\rp $
denote the solution of the following Backward Doubly Stochastic Differential Equation 
(BDSDE) introduced by E.Pardoux and S.Peng in \cite{PPspde}:
\begin{align}
X_t
=&
x+\int_{0}^tb\lp X_s\rp ds
+\int_{0}^t\sigma\lp X_s\rp dW_s,\label{equationX}
\\
 Y_t=&\Phi\lp X_T\rp +\int_t^Tf\XYZs ds
+\int_t^T g\XYs\dBs
-\int_t^TZ_sdW_s.\label{equationYZ}
\end{align}
\subsection*{Assumption}
{\it We suppose that the coefficients $f$ and $g$ satisfy the following:
\begin{align}
\Phi\lp X_T\rp\in&L^2,\nonumber\\
\lv f(x,y,z)-f(x',y',z')\rv^2\leq &L_f\lp |x-x'|^2+|y-y'|^2+|z-z'|^2\rp,\label{cf}\\
\lv g(x,y)-g(x',y')\rv^2\leq &L_g\lp |x-x'|^2+|y-y'|^2\rp,\label{cg}
\end{align}
}
Note that \eqref{cf} and \eqref{cg} yield that $f$ and $g$ have linear
growth in their arguments.
We use two approximations. We at first discretize in time with a constant time mesh $h=T/N$,
 which yields the processes $\lp X^N, Y^N, Z^N\rp$. We then approximate the pair $\lp Y^N, Z^N\rp$ 
by some kind of Picard iteration scheme with $I$ steps $\lp Y^{N,i,I}, Z^{N,I}\rp$ for $i=1,\dots,I$.

In order to be as clear as possible, we introduce below all the definitions used in the paper. 
Most of them are same as in \cite{Gobet}.
\begin{itemize}
\item[{\bf (N0)}] 
{\it For $0\leq t\leq t'\leq T$, set 
$\F_t=
\F^W_t\vee \F^B_{t,T}$ and
\begin{align*}
\F^W_{t}=&\sigma \lp W_s ; 0\leq s\leq t\rp\vee\mathcal N,\quad
\F^B_{t,t'}= \sigma \lp B_s-B_{t'}; t\leq s\leq t'\rp\vee \mathcal N.
\end{align*}
$\E\k$ is the conditionnal expectation with respect to
$\F\tk$.}
\item[{\bf (N1)}] {\it $N$ is the number of steps of the time discretization, 
the integer $I$ corresponds to the number of steps of the Picard iteration, 
$h:=T/N$ is the size of the time mesh and for $k=0,1,\dots,N$ we set $t_k:=kh$ and
$\DeltaB{k}=B\tkk-B\tk$, $\Delta W\kk=W\tkk-W\tk$.
Let $\pi={t_0,t_1,\dots,t_N=T}$ denote the corresponding 
subdivision on $[0,T]$.
}
\item[{\bf (N2) }] {\it The function basis for $\XNt{k}$ is defined as follows:
let $a_k<b_k$ be two reals and $(\mathcal X_i^k)_{i=1...L}$ denote a partition of $[a_k,b_k]$;
for $i=1,\dots,L$ set 
\begin{align}
\uik:=&1_{\mathcal X_i^k}\lp \XNt{k}\rp/\sqrt{P\lp\XNt{k}\in \mathcal X_i^k\rp}
\end{align}
}
\item[{\bf (N3)}] {\it The function basis for $N\sim\mathcal N(0,h)$ is defined as follows:
let $a<b$ two reals and $(\mathcal B_i)_{i=1...L}$ denote a partition of $[a,b]$. 
For $i=1,\dots,L$ set
\begin{align}
\vik{N}:=&1_{\mathcal B_i}\lp N\rp/\sqrt{P\lp N\in \mathcal B_i\rp}
\end{align}
}
\item[{\bf(N4)}] {\it For fixed $k=1,\dots,N$,
let $p_k$ denote the following vector whose components belong to 
$L^2\lp \Omega\rp$. It is defined blockwise as follows:
\begin{align*}
&\lp \uiNk\rp_{i_N},
\lp\uiNk\DeltahB{N-1}\rp_{\iN},
\lp\uiNk\viNmk{\DeltaB{N-1}}\DeltahB{N-2}\rp_{\iN,\iNm},\\
&\dots\\
&\lp\uiNk\prod_{j=k+1}^{N-1}v_{i_{j}}\lp\DeltaB{j}\rp\DeltahB{k}
\rp_{\iN,\iNm,\dots,\ikk}
\end{align*}
where $i_N,\dots,i_{k+1}\in\lb 1,\dots,L\rb$.
Note that $p\k$ is $\F\tk$-measurable and 
$
\E p\k p\k^* = Id
$}
\end{itemize}
\section{Approximation result: step 1}
We first consider a time discretization of equations (\ref{equationX}) and (\ref{equationYZ}).
 The forward equation (\ref{equationX}) is approximated using the Euler scheme:  $\XNt{0}=x$
and for $k=0,\dots,N-1$,
\begin{equation}\label{equationEulerscheme}
 \XNt{k+1}=\XNt{k}+hb(\XNt{k})+\sigma(\XNt{k})\DeltaW{k+1}.
\end{equation}
The following result is well know: (see e.g. \cite{Kloeden})
\begin{theorem}\label{lem-euler-X}
There exists a constant C such that for every $N$
\begin{align*}
\max_{k=1,\dots,N}\sup_{t\km\leq r\leq t\k}\E\lv X_r-X^N\tkm\rvd\leq Ch,
\qquad\max_{k=0,\dots,N}\E\lv X^N\tk\rvd=C<\infty.
\end{align*}
\end{theorem}
The following time regularity is proved in \cite{Aman}  (see also Theorem 2.3 in \cite{Aboura}), 
it extends the original result of Zhang \cite{Zhang}.
\begin{lemma}\label{lem-regularity-Y} 
There exists a constant $C$ 
such that
for every integer $N\geq1$,
$s,t\in[0,T]$,
 \begin{align*}
\sum_{k=1}^N \E\int\tkm^{t_k}\lp\lv Z_r-Z\tkm\rvd+\lv Z_r-Z\tk\rvd \rp dr
\leq C h,\quad
\E\lv Y_t-Y_s\rvd\leq C\lv t-s\rv.
\end{align*} 
\end{lemma}
The backward equation (\ref{equationYZ}) is approximated by backward induction as follows: 
\begin{align}
\YNt{N}:=&\Phi(\XNt{N}),\qquad \ZNt{N}:=0,\label{eq-YZN-tN}\\
\ZNt{k}:=&\unh\Ek\lp\YNt{k+1}\DeltaW{k+1}\rp + \frac{1}{h}\DeltaB{k}\Ek\lp\gXYNkk \DeltaW{k+1}\rp,\label{equationZNL}\\
\YNt{k}:=&\Ek\YNt{k+1}+hf\lp \XNt{k},\YNt{k},\ZNt{k}\rp+\DeltaB{k}\Ek\gXYNkk , \label{equationYN}
\end{align}
Note that as in \cite{Aman}, \cite{Aman2} and \cite{Gobet} we have introduced 
an implicit scheme, thus different from that in
\cite{Aboura}.
However, it differs from that in \cite{Aman}
and \cite{Aman2} since the conditional expectation
we use is taken with respect to $\F\tk$
which is different from $\sigma\lp X^N_{t_j}, j\leq k\rp
\vee\sigma\lp B_{t_j}, j\leq k\rp$ used in \cite{Aman2}.
\begin{proposition}[Existence of the scheme]\label{prop-existance-1}
For sufficiently large N, the above scheme has a unique solution. Moreover, for all $k=0,\dots,N$, we have
$\YNt{k},\ZNt{k}\in L^2\lp\Ftk\rp$.
\end{proposition}
The following theorem is the main result of this section.
\begin{theorem}\label{theorem 1 Gobet}
There exists a constant $C>0$ such that for $h$ small enough
$$
\max_{0\leq k\leq N}\E\lv Y_{t_k}-Y^N\tk\rv^2+\sum_{k=0}^{N-1}\int_{t_k}^{t_{k+1}}\E\lv Z_r-Z^N\tk\rv^2  dr
\leq
 C h+C\E\lv \phi\lp\XNt{N}\rp-\phi\lp X_T\rp\rv^2.
$$
\end{theorem}
The rest of this section is devoted to the proof of this theorem;
it requires several steps.\\
First of all, we define a process $\lp\pY_t,\pZ_t\rp_{t\in[0,T]}$ 
such that $\pY\tk$ and $\pZ\tk$ are $\F\tk$ measurable,
and a family of $\F\tk$ measurable random variables $\uZ\tk$, $k=0,\dots,N$ as follows.
For $t=T$, set 
\begin{align}
\pY_T:=\Phi\lp X^N_{t_N}\rp,\; \pZ_T:=0,\;\uZ_{t_N}:= 0.\label{eq-def-bZ}
\end{align}
Suppose that the scheme $\lp\pY_t,\pZ_t\rp$ is defined for all $t\in [ t_k,T]$ 
and that $\uZ_{t_j}$ has been defined for $j=N,\dots,k$.
Then for $h$ small enough the following equation
\begin{align}
M^k\tkm:=&\E_{k-1}\lp \pY\tk+f\lp X^N\tkm,M^k\tkm,Z^N\tkm\rp\Delta t_{k-1}+
\gXYpk\DeltaB{k-1}\rp \label{eq-M}
\end{align}
has a unique solution.\\
Using Proposition \ref{prop-existance-1} and 
the linear growth of $f$, we deduce that 
the map $F_{\xi}$ defined by
\begin{equation}\label{eq-F-xi}
F_{\xi}(Y)=\E_{k-1} \lp\xi+hf\lp \XNt{k-1},Y,\ZNt{k-1}\rp\rp
\end{equation}
is such that 
$F_{\xi}\lp L^2\lp\F\tkm\rp\rp\subset L^2\lp \F\tkm\rp$.
Futhermore, given $Y,Y'\in L^2\lp\F\tkm\rp$, the $L^2$ contraction property of $\E_{k-1}$
and the Lipschitz condition \eqref{cf} imply
$ \E\lv F_{\xi}(Y)-F_{\xi}\lp Y'\rp\rv^2
\leq h^2L_f\E\lv  Y-Y'\rv^2.
$
Then $F_{\xi}$ is  a contraction for $h$ small enough and 
the fixed point theorem concludes the proof.\\
We can extend $M^k_.$ to the interval $t\in[ t_{k-1},t_k]$ letting
$$
M^k_t:=\E\lp \pY\tk+f\lp X^N\tkm,M^k\tkm,Z^N\tkm\rp\Delta t_{k-1}
+
\gXYpk
\DeltaB{k-1}\right|\left.\F^W_t
\vee\F^B_{t_{k-1},T} \rp,
$$
which is consistent at time $t_{k-1}$. \\
By an extension of the martingale representation theorem
(see e.g. \cite{PPspde} p.212), 
there exists a $\lp\F^W_t\vee\F^B_{t_{k-1},T}\rp_{t_{k-1}\leq t\leq t_k}$-adapted and square integrable process
$\lp N^k_t\rp_{t\in [t_{k-1},t_k]}$ such that for any $t\in [t_{k-1},t_k]$,
$
M^k_t=M^k\tkm+\int\tkm^{t}N^k_sdW_s
$
and hence
$
M^k_t=M^k\tk-\int_t^{t_k}N^k_sdW_s
$. 
Since,
$$
M^k\tk= \pY\tk+f\lp X^N\tkm,M^k\tkm,Z^N\tkm\rp\Delta t_{k-1}+
\gXYpk\DeltaB{k-1},
$$
we deduce that for $t\in[t_{k-1},t_k]$
\begin{equation}\label{eq-MN}
M^k_t= \pY\tk+f\lp X^N\tkm,M^k\tkm,Z^N\tkm\rp\Delta t_{k-1}+
\gXYpk
\DeltaB{k-1}-\int_t^{t_k}N^k_sdW_s.
\end{equation}
For $t\in[ t_{k-1},t_k)$, we set
\begin{equation}\label{eq-def-YZpit}
\pY_t:=M^k_t,\;\pZ_t:=N^k_t,\; \uZ\tkm:=\frac{1}{h}\E_{k-1}\lp\int\tkm^{t_k}\pZ_rdr\rp.
\end{equation}
\begin{lemma}\label{lem-2-6}
 For all $k=0,\dots,N$, 
\begin{equation}\label{eq-YZpN}
 \pY\tk=Y^N\tk,\; \uZ\tk=Z^N\tk
\end{equation}
and hence for $k=1,\dots,N$
\begin{equation}\label{eqYpi}
\pY\tkm= \pY\tk+\int\tkm^{t_k}f\lp X^N\tkm,\pY\tkm,\uZ\tkm\rp dr
+\int\tkm^{t_k}\gXYpk \dBr-\int\tkm^{t_k}\pZ_rdW_r
\end{equation}

\end{lemma}
\proof We proceed  by backward induction. For $k=N$, \eqref{eq-YZpN} is true by \eqref{eq-YZN-tN} and \eqref{eq-def-bZ}. 
Suppose that \eqref{eq-YZpN} holds for  $l=N,N-1,\dots,k$, 
so that $ \pY\tk=Y^N\tk,\; \uZ\tk=Z^N\tk $. 
Then \eqref{eq-YZpN} holds for $l=k-1$; indeed,
for 
$\xi:= \YNt{k}+\DeltaB{k-1} 
\gXYNk
$
we deduce from \eqref{equationYN} and \eqref{eq-M} that 
$M^k\tkm=F_{\xi}\lp M^k\tkm\rp$,
$Y^N\tkm=F_{\xi}\lp Y^N\tkm\rp$ and
$\pY\tkm=M^k\tkm=F_{\xi}\lp M^k\tkm\rp$,
where $F_{\xi}$ is defined by \eqref{eq-F-xi}.
So using the uniqueness of the fixed point of the map $F_{\xi}$,
we can conclude that $\pY\tkm=Y^N\tkm(=M^k\tkm)$.
Therefore, \eqref{eq-MN} and \eqref{eq-def-YZpit} imply \eqref{eqYpi}.
Ito's formula yields 
\begin{align*}
\Delta W\k \int\tkm^{t\k}\pZ_rdW_r=
\int\tkm^{t\k}(W_r-W\tkm)\pZ_rdW_r 
+\int\tkm^{t\k}\int\tkm^r\pZ_sdW_s dW_r 
+\int\tkm^{t\k}\pZ_rdr,
\end{align*}
so that $
\E_{k-1}\lp\Delta W\k\int\tkm^{t_k}\pZ_rdW_r\rp=
\E_{k-1}\lp\int\tkm^{t_k}\pZ_rdr\rp=h\uZ\tkm
$.
Hence multiplying \eqref{eqYpi} by $\Delta W\k$ and taking 
  conditional expectation with respect to $\F\tkm=\F^W\tkm\vee \F^B_{t_{k-1},T}$. We deduce
\begin{align*}
h\uZ\tkm
=&\Ekm\lp Y^N\tk \DeltaW{k} \rp
+\DeltaB{k-1}\Ekm\lp 
\gXYNk
 \DeltaW{k}\rp
\end{align*}
Comparing this with \eqref{equationZNL} concludes the proof of \eqref{eq-YZpN} for $l=k-1$.
\endproof
Lemma \ref{lem-2-6} shows that for $r\in[t\k,t\kk]$ one can upper estimate the 
$L^2$ norm of $Z_r-Z^N\tk$ by that of
$Z_r-\pZ_r$ and increments of $Z$.
Indeed, using \eqref{eq-YZpN} we have for $k=0,\dots,N-1$ and $r\in[t_k,t\kk]$
\begin{align*}
\E\lv Z_r-Z^N\tk\rvd=\E\lv Z_r-\uZ\tk\rvd
\leq 2\E\lv Z_r-Z\tk\rvd+2\E\lv Z\tk-\uZ\tk\rvd
\end{align*}
Furthermore, \eqref{eq-def-YZpit} and Cauchy-Schwarz's inequality yield for $k=0,\dots,N-1$
\begin{align*}
\E\lv Z\tk-\uZ\tk\rvd
\leq & \unh\E\int\tk^{t\kk}\lv Z\tk-\pZ_r\rvd dr\\
\leq & \frac{2}{h} \E\int\tk^{t\kk}\lv Z\tk-Z_r\rvd dr+\frac{2}{h} \E\int\tk^{t\kk}\lv Z_r-\pZ_r\rvd dr.
\end{align*}
Hence we deduce
\begin{align}
\sum_{k=0}^{N-1}\int\tk^{t\kk}\E\lv Z_r-Z^N\tk\rvd dr
\leq &6\sum_{k=0}^{N-1}\int\tk^{t\kk}\E\lv Z_r-Z\tk\rvd dr
+4\sum_{k=0}^{N-1}\int\tk^{t\kk}\E\lv Z_r-\pZ_r\rvd dr. \label{eq-numero2}
\end{align}
Using Lemma \ref{lem-regularity-Y} and \eqref{eq-numero2} we see that
Theorem \ref{theorem 1 Gobet} is a straightforward consequence of the following:
\begin{theorem}\label{theorem 1 Gobet bis}
There exists a constant C such that for $h$ small enough,
\begin{align*}
\max_{0\leq k\leq N}\E\lv Y_{t_k}-\pY\tk\rv^2+\int_{0}^{T}\E\lv Z_r-\pZ_r\rv^2  dr
\leq
 C h+C\E\lv \Phi\lp\XNt{N}\rp-\Phi\lp X_T\rp\rv^2.
\end{align*}
\end{theorem}
\proof
For any $k=1,\dots,N$ set 
\begin{equation}\label{Ikm}
I_{k-1}:=\E\lv Y\tkm-\pY\tkm\rv^2+\E\int_{t_{k-1}}^{t_k}\lv Z_r-\pZ_r\rv^2 dr.
\end{equation}

Since $Y\tkm-\pY\tkm$ is $\F\tkm$-measurable while 
for $r\in[t\k,t\kk]$ the random variable $ Z_r-\pZ_r$ 
is $\F^W_r\vee\F^B_{t\km,T}$-measurable, we deduce that
 $Y\tkm-\pY\tkm$ is orthogonal to $\int_{t_{k-1}}^{t_k}\lp Z_r-\pZ_r\rp dW_r$.
Therefore,
the identities \eqref{equationYZ} and \eqref{eqYpi} imply that
\begin{align*}
 I_{k-1}=&\E\lv Y\tkm-\pY\tkm+\int_{t_{k-1}}^{t_k}\lp Z_r-\pZ_r\rp dW_r\rv^2\\
=&\E\lv Y\tk-\pY\tk+\int_{t_{k-1}}^{t_k}\lp f\XYZr-f\XYZpikm \rp dr\right.\\
&\qquad\qquad\left.+\int_{t_{k-1}}^{t_k} \lp \gXYr  -\gXYpk\rp\dBr\rv^2.
\end{align*}
Notice that for $t\km\leq r\leq t\k$ the random variable 
$ \gXYr - \gXYpk$ is $\F\tk^W\vee\F^B_{r,T}$-measurable.
Hence 
$Y\tk -\pY\tk$, which is $\F\tk$-measurable,
and $ \int_{t_{k-1}}^{t_k} \lp \gXYr - \gXYpk\rp\dBr$
are orthogonal. The inequality
$
(a+b+c)^2\leq \lp 1+\frac{1}{\lambda}\rp (a^2+c^2)+\lp 1+2\lambda\rp b^2 +2ac
$
valid for $\lambda>0$,
Cauchy-Schwarz's inequality and the isometry of backward stochastic integrals yield
for $\lambda:=\frac{\e}{h}$, $\e>0$:
\begin{align*}
 I_{k-1}
\leq &
\lp 1+\frac{h}{\e}\rp\lc\E\lv Y\tk-\pY\tk\rv^2+\E\lv\int_{t_{k-1}}^{t_k}\lp \gXYr - \gXYpk\rp\dBr\rv^2\rc\\
&+\lp1+2\frac{\e}{h}\rp\E\lv\int_{t_{k-1}}^{t_k}\lp f\XYZr-f\XYZpikm\rp dr\rv^2\\
\leq & 
\lp 1+\frac{h}{\e}\rp
\lc\E\lv Y_k-\pY\tk\rv^2+\E\int_{t_{k-1}}^{t_k}\lv \gXYr - \gXYpk \rv^2dr\rc\\
&+\lp h+2\e\rp\E\int_{t_{k-1}}^{t_k}\lv f\XYZr-f\XYZpikm\rv^2 dr.
\end{align*}
The Lipschitz properties \eqref{cf} and \eqref{cg} of $f$ and $g$ imply
\begin{align}
I_{k-1}
\leq & 
\lp 1+\frac{h}{\e}\rp \lc\E\lv Y\tk-\pY\tk\rv^2+
 L_g \E\int_{t_{k-1}}^{t_k}\lp\lv X_r-X^N\tk\rv^2+\lv Y_r-\pY\tk \rv^2\rp dr\rc\nonumber\\
&+\lp h+2\e\rp L_f\E\int_{t_{k-1}}^{t_k}\lp\lv X_r-X^N\tkm\rv^2+\lv Y_r-\pY\tkm \rv^2+\lv Z_r-\uZ\tkm \rv^2\rp dr.
\label{star-star}
\end{align}
Using the definition of $\uZ\tk$ in \eqref{eq-def-YZpit}, the $L^2$ contraction property of $\Ek$ 
and Cauchy-Schwarz's inequality, we have
\begin{align*}
h\E\lv Z\tk-\uZ\tk\rv^2
\leq  \frac{1}{h}\E\lv\Ek\lp\int_{t_k}^{t_{k+1}}\lp Z\tk-\pZ_r\rp dr\rp\rv^2
\leq  \E\int_{t_k}^{t_{k+1}}\lv Z\tk-\pZ_r\rv^2 dr.
\end{align*}
Thus, by Young's inequality,
we deduce for 
$k=1,\dots,N$
\begin{align*}
\E\int_{t_{k-1}}^{t_k}\lv Z_r-\uZ\tkm \rv^2 dr
\leq &
2\E\int_{t_{k-1}}^{t_k}\lv Z_r-Z\tkm \rv^2 dr
+2h\E\lv Z\tkm-\uZ\tkm \rv^2 \\
\leq &
2\E\int_{t_{k-1}}^{t_k}\lv Z_r-Z\tkm \rv^2 dr
+4\E\int\tkm^{t_{k}}\lv\pZ_r-Z_r\rv^2dr\\
&+4\E\int\tkm^{t_{k}}\lv Z_r-Z\tkm\rv^2dr\\
\leq &
6\E\int_{t_{k-1}}^{t_k}\lv Z_r-Z\tkm \rv^2 dr
+4\E\int\tkm^{t_{k}}\lv\pZ_r-Z_r\rv^2dr.
\end{align*}
We now deal with increments of $Y$. Using Lemma \ref{lem-regularity-Y}, we have
\begin{align*}
\E\int_{t_{k-1}}^{t_k}\lv Y_r-\pY\tkm \rv^2 dr\leq &
2\E\int_{t_{k-1}}^{t_k}\lv Y_r-Y\tkm \rv^2 dr 
+2\E\int_{t_{k-1}}^{t_k}\lv Y\tkm-\pY\tkm \rv^2 dr  \\
\leq &
Ch^2
+2h\E\lv Y\tkm-\pY\tkm \rv^2,
\end{align*}
while a similar argument yields
\begin{align*}
\E\int_{t_{k-1}}^{t_k}\lv Y_r-\pY\tk \rv^2 dr \leq &
Ch^2
+2h\E\lv Y\tk-\pY\tk \rv^2.
\end{align*}
Using Theorem \ref{lem-euler-X} and the previous upper estimates in \eqref{star-star}, we deduce
\begin{align*}
I_{k-1}
\leq &
 \lp 1+\frac{h}{\e}\rp\E\lv Y\tk-\pY\tk\rv^2
+L_f\lp h+2\e\rp
\lc
Ch^2+2 h\E\lv Y\tkm-\pY\tkm \rv^2\right. \\
&+\left.6
\E\int_{t_{k-1}}^{t_k}\lv Z_r-Z\tkm \rv^2 dr
+4\E\int\tkm^{t_{k}}\lv\pZ_r-Z_r\rv^2dr\rc\\
&+L_g\lp 1+\frac{h}{\e}\rp  
\lc
Ch^2
+2  h\E\lv Y\tk-\pY\tk \rv^2\rc.
\end{align*}
Thus, \eqref{Ikm} implies that for any $\e>0$
\begin{align*}
&\lc 1-2L_f\lp h+2\e\rp h\rc \E\lv Y\tkm-\pY\tkm\rv^2
+\lc1-4L_f\lp h+2\e\rp \rc\E\int_{t_{k-1}}^{t_k}\lv Z_r-\pZ_r\rv^2 dr\\
\leq & 
\lp 1+\frac{h}{\e}+2L_g\lp 1+\frac{h}{\e}\rp  h\rp\E\lv Y\tk-\pY\tk\rv^2
+\lp L_f\lp h+2\e\rp +L_g\lp 1+\frac{h}{\e}\rp \rp Ch^2 \\
&+6L_f\lp h+2\e\rp\E\int\tkm^{t_{k}}\lv Z_r-Z\tkm\rv^2dr.
\end{align*}
Now we choose $\e$ such that $8\e L_f=\undeux$. 
Then we have for $\tC=4L_f$, $h$ small enough and some positive constant $\bC$ 
depending on $L_f$ and $L_g$:
\begin{align}
\lp 1-\tC h\rp &\E\lv Y\tkm-\pY\tkm\rv^2
+\lp\undeux-\tC h \rp\E\int_{t_{k-1}}^{t_k}\lv Z_r-\pZ_r\rv^2 dr\nonumber\\
\leq & 
\lp 1+\bC h\rp\E\lv Y\tk-\pY\tk\rv^2   + \bC h^2 
 +\bC\E\int\tkm^{t_{k}}\lv Z_r-Z\tkm\rv^2dr.\label{Y-Ypi}
\end{align}
We need the following
\begin{lemma}\label{lem-1-Ch}
Let $L>0$; then for $h^*$ small enough (more precisely $Lh^*<1$) 
there exists $\Gamma:=\frac{L}{1-Lh^*}>0$ such that for all $h\in (0,h^*)$
we have
$
\frac{1}{1-Lh}< 1+\Gamma h
$
\end{lemma}
\proof
Let  $h\in(0,h^*)$; then we have $1-Lh>1-Lh^*>0$. 
Hence $ \frac{L}{1-Lh}<\frac{L}{1-Lh^*}=\Gamma$, so that 
$Lh<\Gamma h(1-Lh)$,
which yields $ 1+\Gamma h-L h-\Gamma Lh^2=(1+\Gamma h)(1-Lh)>1$.
This concludes the proof.
\endproof
Lemma \ref{lem-1-Ch} and \eqref{Y-Ypi} imply the existence of a constant $C>0$ such that
for $h$ small enough and $k=1,2,\dots,N$ we have
\begin{align}
\E\lv Y\tkm-\pY\tkm\rv^2
\leq  
\lp 1+C h\rp\E\lv Y\tk-\pY\tk\rv^2   + Ch^2 
+C\E\int\tkm^{t_{k}}\lv Z_r-Z\tkm\rv^2dr.
\end{align}
The final step relies on the following discrete version of Gronwall's lemma (see \cite{Gobet}).
\begin{lemma}[Gronwall's Lemma]\label{lem-Gronwall}
Let $(a_k),(b_k),(c_k)$  be nonnegative sequences such that for some $K>0$ we have for all $k=1,\dots,N-1$,
$a_{k-1}+c_{k-1}\leq (1+Kh)a_k+b_{k-1}$. Then, for all $k=0,\dots,N-1$,
$
a_k+\sum_{i=k}^{N-1}c_i\leq e^{K(T-t_k)}\lp a_N+\sum_{i=k}^{N-1}b_i\rp
$
\end{lemma}
Use Lemma \ref{lem-Gronwall} with $c\k=0$,  $a\km=\E\lv Y\tkm-\pY\tkm\rvd$ and 
$b\k=C\E\int\tkm^{t\k}\lv Z_r-Z\tkm\rvd+Ch^2$;
this yields
\begin{align}
\sup_{0\leq k\leq N}\E \lv Y\tk-\pY\tk\rv^2\leq &C\lp\E \lv Y_T-\pY_{t_N}\rv^2
+\sum_{k=1}^{N}\E\int\tkm^{t_{k}}\lv Z_r-Z\tkm\rv^2dr +Ch\rp\nonumber\\
\leq & C\lp\E \lv Y_T-\pY_{t_N}\rv^2+Ch\rp,\label{majoY}
\end{align}
where the last upper estimate is deduced from Lemma \ref{lem-regularity-Y}.
We sum \eqref{Y-Ypi} from $k=1$ to $k=N$; using \eqref{majoY} we deduce 
that for some constant $\bar C$ depending on $L_f$ and $L_g$ we have
\begin{align*}
\lp\undeux-\tC h \rp\E\int_{0}^{T}\lv Z_r-\pZ_r\rv^2 dr
\leq & 
 \bC h\sum_{k=1}^{N-1}\E\lv Y\tk-\pY\tk\rv^2   + \bC h + \bC \E \lv Y_T-\pY_{t_N}\rv^2\\
\leq & \bC h + \bC \E \lv Y_T-\pY_{t_N}\rv^2+\bC h \lp \bC+N\E \lv Y_T-\pY_{t_N}\rv^2 \rp\\
\leq & \bC h + \bC\E \lv Y_T-\pY_{t_N}\rv^2.
\end{align*}
The definitions of $Y_T$ and $Y^N_{t_N}$ from \eqref{equationYZ} and \eqref{eq-YZN-tN} conclude the proof of Theorem \ref{theorem 1 Gobet bis}.
\endproof
\section{Approximation results: step 2}
In order to approximate $\lp \YNt{k},\ZNt{k}\rp_{k=0,\dots,N}$ we use the idea of E.Gobet,
J.P. Lemor and X.Warin \cite{Gobet},
that is a projection on the function basis and a Picard iteration scheme. In this section, $N$ and $I$ are fixed positive integers.
We define the sequences $\lp\YNIt{i}{k}\rp_{i=0,...,I\; k=0,...,N}$ and $\lp Z^{N,I}\tk\rp_{k=0,...,N-1}$  
using backward induction on $k$, and for fixed $k$ forward induction on $i$ for $Y^{N,i,I}\tk$ as follows:
For $k=N$, $Z\NI_{t_N}=0$ and for $i=0,\dots,I,$ set
$\YNIt{i}{N}:=P_N\Phi\lp\XNt{N}\rp$.
Assume that $Y\NII\tkk$ has been defined and set 
 \begin{align}
Z^{N,I}\tk:=& \frac{1}{h} P_{k} \lc Y^{N,I,I}\tkk\Delta W_{k+1}\rc
+\frac{1}{h} P_{k}\lc\DeltaB{k}\gTNIkk\DeltaW{k+1}\rc\label{eq-ZNI}.
\end{align}
Let $\YNIt{0}{k}:=0$ and for $i=1,\dots,I$ define inductively
by the following Picard iteration scheme:
\begin{align}
 Y\NiI\tk:=& P_{k}Y^{N,I,I}\tkk +hP_k\lc f\lp X^N_{t_{k}},Y^{N,i-1,I}\tk,Z^{N,I}\tk\rp\rc+P_{k}\lc\DeltaB{k}\gTNIkk\rc\label{eq-YNiI},
 \end{align}
where $P_k$ is the orthogonal projection on the Hilbert space $\P_k\subset L^2\lp\F\tk\rp$ 
generated by the function $p_k$ defined by (N4).
Set $R_k:=I-P_k$. Note that $P\k$ is a contraction of $L^2\lp\F\tk\rp$. Furthermore, given $Y\in L^2\lp\Omega\rp$,
\begin{equation}\label{lem-PE=P}
\Ek P\k Y=P\k\E\k Y=P\k Y.
\end{equation}
Indeed, since $\P\k\subset L^2\lp \F\tk\rp$,
$\Ek P\k Y=P\k Y$.
Let $Y\in L^2$;
for every, $U_k\in \P_k$, since $U_k$ is $\F\tk$-measurable,
 we have $\E\lp U_kR_kY\rp=0=\E\lp U_k\Ek R_kY\rp$; so that, $P_k\Ek R_k(Y)=0$.
Futhermore $Y=P_kY+R_kY$ implies $P_k\Ek Y=P_kP_kY+P_k\Ek R_k Y=P_kY$ which yields \eqref{lem-PE=P}.
Now we state the main result of this section.
\begin{theorem}\label{th-2}For $h$ small enough, we have
\begin{align*}
\max_{0\leq k\leq N}&\E\lv\YNIt{I}{k}-\YNt{k}\rv^2
+h\sum_{k=0}^{N-1}\E\lv Z\NI\tk-\ZNt{k}\rv^2
\leq C h^{2I-2}+C\sum_{k=0}^{N-1}\E\lv R_k Y^N\tk\rv^2\\
& +C\E \lv \Phi\lp X^N_{t_N}\rp - P_N\Phi\lp X^N_{t_N}\rp \rvd 
+Ch\sum_{k=0}^{N-1}\E\lv R_k Z^{N,I}\tk\rv^2.
\end{align*}
\end{theorem}
\subsection*{Proof of Theorem \ref{th-2}}
The proof will be deduced from severals lemmas. 
The first result gives integrability properties of the scheme defined by \eqref{eq-ZNI} and \eqref{eq-YNiI}.
\begin{lemma}For every $k=0,\dots,N$ and $i=0,\dots,I$ we have
$Y^{N,i,I}_{t_k},Z^{N,I}_{t_k}\in L^2\lp \F_{t_k}\rp$.
\end{lemma}
\proof We prove this by backward induction on $k$, and for fixed $k$ by forward induction on $i$.
By definition $Y^{N,i,I}_{t_N}=P_N\Phi(X^N_{t_N})$ and $Z^{N,I}_{t_N}=0$. 
Suppose that $Z\NI_{t_j}$ and $Y^{N,l,I}_{t_j}$ belong to  $L^2\lp \F_{t_j}\rp$ for $j=N,N-1,\dots,k+1$
and any $l$,
and for $j=k$ and $l=0,\dots,i-1$;
we will show that $Y^{N,i,I}_{t_k},Z^{N,I}_{t_k}\in L^2\lp \F_{t_k}\rp$.\\
The measurability is obvious since $\mathcal P\k\subset L^2\lp \F\k\rp$.
We at first prove the square integrability of $Z^{N,I}_{t_k}$.
Using \eqref{lem-PE=P}, the conditional Cauchy-Schwarz inequality and the independence of 
$\Delta W\kk$ and $\F\tk$, we deduce
\begin{align*}
\E\lv P_k\lp Y\NII\tkk\Delta W_{k+1}\rp\rv^2
=&
\E\lv P_k\Ek\lp Y\NII\tkk\Delta W_{k+1}\rp\rv^2
\leq 
\E\lv \Ek\lp Y\NII\tkk\Delta W_{k+1}\rp\rv^2\\
\leq &
\E\lp \Ek\lv\Delta W\kk\rvd\Ek\lv Y\NII\tkk\rv^2\rp
\leq 
h\E \lv Y\NII\tkk\rv^2.
\end{align*}
A similar computation using the independence of $\Delta W\kk$ and $\F\tk$,
and of $\DeltaB{k}$ and $\F\tkk$ as well as the growth condition 
deduced from \eqref{cg} yields
\begin{align*}
\E&\lv P_k\lp \DeltaB{k}\Delta W_{k+1} \gTNIkk\rp\rv^2
= 
\E\lv P_k\Ek\lp \DeltaB{k}\Delta W_{k+1}\gTNIkk\rp\rv^2\\
\leq &
\E\lv \Ek\lp \DeltaB{k}\Delta W_{k+1}\gTNIkk\rp\rv^2
\leq 
h\E\E\kk\lv \DeltaB{k}\gTNIkk\rv^2\\
\leq &
h^2\E \lv\gTNIkk\rv^2
\leq 
2h^2 \lv g(0,0)\rv^2 +2h^2L_g\lp\E\lv X^N\tkk\rv^2+\E\lv Y\NII\tkk\rv^2\rp.
\end{align*}
The two previous upper estimates and the induction hypothesis proves that $Z\NI\tk\in L^2\lp\F\tk\rp$.
A similar easier proof shows that $Y\NiI\tk\in L^2\lp\F\tk\rp$.
\endproof
The following lemma gives $L^2$ bounds for multiplication by $\Delta W\kk$
\begin{lemma}\label{lem-Y-DeltaW}
For every $Y\in L^2$ we have
$
\E\lv \E_k\lp Y\Delta W_{k+1}\rp\rv^2
\leq  h \lp\E |Y|^2-\E\lv\E_kY\rv^2\rp
$
\end{lemma}
\proof
Using the fact that $\E_k\lp \Delta W_{k+1}\E_k Y\rp=0$ we have
\begin{align*}
\E\lv \E_k\lp Y\Delta W_{k+1}\rp\rv^2
=&\E\lv \E_k\lp\lp Y-\E_kY\rp \Delta W_{k+1}\rp\rv^2
\end{align*}
Using the conditional Cauchy-Schwarz inequality and the independence of $\Delta W\kk$ 
and $\F\tk$, we deduce
$
\E\lv \E_k\lp Y\Delta W_{k+1}\rp\rv^2
\leq
 h\E \lv Y-\E_kY\rv^2
\leq 
 h \lp\E |Y|^2-\E\lv\E_kY\rv^2\rp;
$
this concludes the proof.
\endproof
The following result gives orthogonality properties
of several projections.
\begin{lemma}\label{lem-1}
 Let $k=0,\dots,N-1$, and $M\tkk,N\tkk\in L^2\lp\F\tkk\rp$. Then
$$ 
\E\lp P_k M\tkk P_k\lp\DeltaB{k}N\tkk\rp\rp=0.
$$
\end{lemma}
\proof
Let $M\tkk\in L^2\lp\F\tkk\rp$; the definition of $P\k$ yields
\begin{align}
&P\k M\tkk=\sum_{1\leq i_N\leq L}\alpha(i_N)u_{i_N}\lp X^N\tk\rp
+\sum_{1\leq i_N\leq L}\alpha(N-1,i_N)u_{i_N}\lp X^N\tk\rp\DeltahB{N-1}\label{dec-1}\\
&+\sum_{k\leq l\leq N-1}\sum_{1\leq i_N,\dots,i_{l+1}\leq L}
\alpha\lp l,i_N,\dots,i_{l+1}\rp u_{i_N}\lp X^N\tk\rp 
\prod_{r=l+1}^{N-1} 
v_{i_{r}}\lp \DeltaB{r}\rp
\DeltahB{l},\nonumber
\end{align}
where $\alpha(i_N)=\E\lc M\tkk u_{i_N}\lp X^N\tk\rp\rc$, 
$\alpha\lp N-1,i_N\rp=\E\lc M\tkk u_{i_N}\lp X^N\tk\rp\DeltahB{N-1}\rc$,
and\\ 
$\alpha\lp l,i_N,\dots,i_{l+1}\rp=\E\lc M\tkk u_{i_N}\lp X^N\tk\rp
\prod_{r=l+1}^{N-1} v_{i_{r}}\lp\DeltaB{r}\rp
\DeltahB{l}\rc$.
Taking conditional expectation with respect to $\F\tkk$, we deduce that for any 
$i_N,\dots,i_k\in\lbrace 1,\dots,L\rbrace$
\begin{align*}
\alpha\lp k,i_N,\dots,i\kk\rp=&
\E\lc M\tkk \uiNk\prod_{r=k+1}^{N-1} 
v_{i_{r}}\lp\DeltaB{r}\rp\E_{k+1}\DeltahB{k}\rc=0.
\end{align*}
A similar decomposition of $P\k\lp\DeltaB{k}N\tkk\rp$ yields
\begin{align}
& P\k\lp\DeltaB{k}N\tkk\rp
=\sum_{1\leq\iN\leq L}\beta\lp i_N\rp \uiNk
+\sum_{1\leq\iN\leq L}\beta\lp N-1,i_N\rp \uiNk\DeltahB{N-1}\nonumber
\\
&+\sum_{k\leq l\leq N-1}\sum_{1\leq i_N,\dots,i_{l+1}\leq L}
\beta\lp l,i_N,\dots,i_{l+1} \rp \uiNk
\prod_{r=l+1}^{N-1}
v_{i_{r}}\lp\DeltaB{r}\rp\DeltahB{l}\label{dec-2}
\end{align}
where $\beta\lp i_N\rp=\E\lc \DeltaB{k}N\tkk u_{i_N}\lp X^N\tk\rp\rc$,
$\beta\lp N-1,i_N\rp=\E\lc \DeltaB{k}N\tkk  u_{i_N}\lp X^N\tk\rp\DeltahB{N-1}\rc$ 
and $\beta\lp l,i_N,\dots,i_{l+1}\rp=\E\lc \DeltaB{k}N\tkk  u_{i_N}\lp X^N\tk\rp
\prod_{r=l+1}^{N-1}
v_{i_{r}}\lp \DeltaB{r}\rp\DeltahB{l}\rc$.
In the above sum, all terms except those corresponding to $l=k$ are equal to 0.
Indeed, let $l\in\lbrace k+1,\dots,N-1\rbrace$; then using again the 
conditional expectation with respect to $\F\tkk$
we obtain
\begin{align*}
\beta\lp l,i_N,\dots,i_{l+1} \rp=
&\E\lc N\tkk \uiNk\prod_{r=l+1}^{N-1}
v_{i_{r}}\lp\DeltaB{r}\rp\DeltahB{l}\E\kk\DeltaB{k}\rc
=0
\end{align*}
The two first terms in the decomposition of $P\k\lp\DeltaB{k}N\tkk\rp$ are dealt with by a similar argument.
Notice that for any $l\in\lbrace k+1,\dots,N-1\rbrace$ and any
$i_N,\dots,i_l,j_N,\dots,j\kk\in\{1,\dots,L\}$ we have,
(conditioning with respect to $\F\tk$):
\begin{align*}
\E\lc \uiNk
\prod_{r=l+1}^{N-1}
v_{i_{r}}\lp\DeltaB{r}\rp\DeltahB{l}
u_{j_N}\lp X^N\tk\rp \prod_{r=k+1}^{N-1}
v_{j_{r}}\lp\DeltaB{r}\rp\DeltahB{k}\rc=0.
\end{align*}
A similar computation proves that for any $i_N,j_N,\dots,j\kk
\in\{1,\dots,L\}$, $\xi\in\left\{1,\DeltahB{N-1}\right\}$
\begin{align*}
\E\lc u_{i_N}\lp X^N\tk\rp \xi u_{j_N}\lp X^N\tk\rp\prod_{r=k+1}
^{N-1}v_r\lp \DeltaB{r}\rp\rc=0.
\end{align*}
The decompositions \eqref{dec-1} and \eqref{dec-2} conclude the proof.
\endproof
The next lemma provides upper bounds of the $L^2$-norm
of $Z\NI\tk$ and $Z\NI\tk-Z^N\tk$.
\begin{lemma}\label{lem-2} For small $h$ enough and for $k=0,\dots,N-1$,
we have the following $L^2$ bounds\\
\begin{align}
\E\lv Z^{N,I}\tk\rv^2\leq& \frac{1}{h}\lp \E\lv  Y^{N,I,I}\tkk \rv^2-\E\lv\Ek  Y^{N,I,I}\tkk \rv^2 \rp\label{bound-ZNI}\\
&+\frac{1}{h}\lp \E\lv \DeltaB{k} \gTNIkk\rv^2 - \E\lv \E_k\lp\DeltaB{k} \gTNIkk\rp\rv^2 \rp\nonumber,
\end{align}
\begin{align}
\E\lv Z^{N,I}\tk-Z^N\tk\rv^2
\leq &
\E\lv R_kZ^N\tk\rv^2+\frac{1}{h}\lp \E\lv \YNIt{I}{k+1}-Y^N\tkk\rv^2 -\E\lv\Ek\lp \YNIt{I}{k+1}-Y^N\tkk\rp\rv^2 \rp\label{bound-ZNI-ZN}\\
&+\frac{1}{h}\lp \E\lv\DeltaB{k}\lc\gTNIkk-\gXYNkk\rc \rv^2\right.\nonumber\\
&\left.\qquad- \E\lv\Ek\lp \DeltaB{k}\lc\gTNIkk-\gXYNkk\rc \rp\rv^2\rp.\nonumber
\end{align}
\end{lemma}
\proof 
Lemma \ref{lem-1} implies that both terms in the right hand side
of \eqref{eq-ZNI} are orthogonal. Hence
squaring both sides of equation \eqref{eq-ZNI}, 
 using \eqref{lem-PE=P} and Lemma \ref{lem-Y-DeltaW}, we deduce
\begin{align*}
 \E\lv Z^{N,I}\tk\rv^2
=&\frac{1}{h^2}\E\lv P_{k}\lp Y^{N,I,I}\tkk\Delta W_{k+1}\rp\rv^2
+\frac{1}{h^2}\E\lv P_k\lp\DeltaB{k} \gTNIkk\Delta W_{k+1}\rp\rv^2\\
=&\frac{1}{h^2}\E\lv P_{k}\Ek\left[  Y^{N,I,I}\tkk\Delta W_{k+1}\right] \rv^2
+\frac{1}{h^2}\E\lv P_k\E_k\lp\DeltaB{k} \gTNIkk \Delta W_{k+1}\rp\rv^2\\
\leq& \frac{1}{h^2}\E\lv \Ek\left[  Y^{N,I,I}\tkk\Delta W_{k+1}\right] \rv^2 
+\frac{1}{h^2}\E\lv \E_k\lp\DeltaB{k} \gTNIkk\Delta W_{k+1}\rp\rv^2 \\
\leq &\frac{1}{h}\lp \E\lv  Y^{N,I,I}\tkk \rv^2-\E\lv\Ek  Y^{N,I,I}\tkk \rv^2 \rp\\&+\frac{1}{h}\lp \E\lv \DeltaB{k} \gTNIkk\rv^2 - \E\lv \E_k\lp\DeltaB{k} \gTNIkk\rp\rv^2 \rp ;
\end{align*}
this proves \eqref{bound-ZNI}.\\
Using the orthogonal decomposition $Z^N\tk=P_kZ^N\tk+R_kZ^N\tk$, 
since $Z\NI\tk\in\mathcal P\k$ we have
$\E\lv Z^{N,I}\tk-Z^N\tk\rv^2=
\E\lv Z^{N,I}\tk-P_kZ^N\tk\rv^2+\E\lv R_kZ^N\tk\rv^2.$
Futhermore \eqref{equationZNL}, \eqref{eq-ZNI} and \eqref{lem-PE=P} yield
\begin{align*}
Z^{N,I}\tk-P_k Z^N\tk=&\frac{1}{h}P_k\lc \lp\YNIt{I}{k+1}-\YNt{k+1}\rp\DeltaW{k+1}\rc\\
&+\frac{1}{h}P_k\lc \DeltaB{k}\lp \gTNIkk-\gXYNkk  \rp\DeltaW{k+1}\rc.
\end{align*}
Lemma \ref{lem-1} shows that the above decomposition is 
orthogonal; thus using \eqref{lem-PE=P}, the contraction property
of $P\k$ and Lemma \ref{lem-Y-DeltaW}, we deduce
\begin{align*}
\E\lv Z^{N,I}\tk-P_k Z^N\tk\rv^2
=&\frac{1}{h^2}\E\lv P_k\Ek\lc \lp\YNIt{I}{k+1}-\YNt{k+1}\rp\DeltaW{k+1}\rc\rv^2\\
&+\frac{1}{h^2}\E\lv P_k\Ek\lc \DeltaB{k}\DeltaW{k+1}\lp \gTNIkk-\gXYNkk   \rp\rc\rv^2\\
\leq&\frac{1}{h^2}\E\lv \Ek\lc \lp\YNIt{I}{k+1}-\YNt{k+1}\rp\DeltaW{k+1}\rc\rv^2\\
&+\frac{1}{h^2}\E\lv \Ek\lc \DeltaB{k}\DeltaW{k+1}\lp \gTNIkk-\gXYNkk  \rp\rc\rv^2\\
\leq &\frac{1}{h}\lp \E\lv \YNIt{I}{k+1}-Y^N\tkk\rv^2 -\E\lv\Ek\lp \YNIt{I}{k+1}-Y^N\tkk\rp\rv^2 \rp\\
&+\frac{1}{h}\lp \E\lv\DeltaB{k}\lc\gTNIkk-\gXYNkk\rc \rv^2\right.\\
&\left.\qquad- \E\lv\Ek\lp \DeltaB{k}\lc\gTNIkk-\gXYNkk\rc \rp\rv^2\rp.
\end{align*}
This concludes the proof of \eqref{bound-ZNI-ZN}.
\endproof
For $Y\in L^2\lp \F\tk\rp$, let 
$\chi_k^{N,I}(Y)$ be defined by:
\begin{align*}
\chi_k^{N,I}(Y):=P_{k}\lp Y^{N,I,I}\tkk+hf\lp X^N\tk,Y,Z^{N,I}\tk\rp + \DeltaB{k}\gTNIkk\rp.
\end{align*}
The growth conditions of $f$ and $g$ deduced from \eqref{cf}, \eqref{cg} and the 
orthogonality of $\DeltaB{k}$ and $\F\tkk$ imply that
$\chi\NI\k\lp L^2\lp\F\tk\rp\rp\subset\mathcal P\k\subset L^2\lp\F\tk\rp$.
Futhermore, \eqref{cf} implies that for $Y_1,Y_2\in L^2\lp\F\tk\rp$
\begin{align}
 \E\lv \chi^{N,I}_k(Y_2)-\chi^{N,I}_k(Y_1) \rv^2
\leq L_fh^2\E\lv Y_2-Y_1 \rv^2, \label{L3.7}
\end{align}
and \eqref{eq-YNiI} shows that $Y^{N,i,I}\tk=\chi^{N,I}\k\lp Y^{N,i-1,I}\tk\rp$
for $i=1,\dots,I$.
\begin{lemma}\label{lem-3-8}
 For small $h$ (i.e., $h^2L_f<1$) and for $k=0,\dots,N-1$,
there exists a unique
$Y^{N,\infty,I}\tk\in L^2\lp\F\tk\rp$ such that 
\begin{align}
&Y^{N,\infty,I}\tk=
P_{k}\lc Y^{N,I,I}\tkk+hf\lp X^N\tk,Y\NooI\tk,Z^{N,I}\tk\rp 
+\DeltaB{k}\gTNIkk\rc,\label{eq-Yinf}\\
&\E\lv Y^{N,\infty,I}\tk-Y^{N,i,I}\tk\rv^2\leq 
L_f^ih^{2i}\E\lv Y^{N,\infty,I}\tk\rv^2,\label{error-Yinf-YNiI}
\end{align}
and there exists some constant $K>0$ such that for every
$N,k,I$,
\begin{equation}\label{bound-Yinf}
 \E\lv Y^{N,\infty,I}\tk\rv^2\leq Kh+(1+Kh)\E\lv Y^{N,I,I}\tkk\rv^2.
\end{equation}
\end{lemma}
\proof
The fixed point theorem applied to the map $\chi\NI\k$,
which is a contration for $h^2L_f<1$, proves \eqref{eq-Yinf} ;
 \eqref{error-Yinf-YNiI} is straightforward consequence from
\eqref{eq-YNiI} by induction on $i$.
Lemma \ref{lem-1} shows that $P\k Y\NII\tkk$ and
$P\k\lp \DeltaB{k} g\lp X^N\tkk,Y\NII\tkk\rp\rp$ are orthogonal.
Hence for any $\e>0$,
using Young's inequality, \eqref{lem-PE=P}, the $L^2$ contracting property of $P\k$, 
the growth condition of $g$ deduced from \eqref{cg} we obtain
\begin{align*}
\E\lv Y^{N,\infty,I}\tk\rv^2
\leq 
&\lp1+\frac{ h}{\e}\rp \E\lv P_{k} Y^{N,I,I}\tkk \rv^2 
+(h^2 +2\e h)\E\lv P_{k}\lc f\lp X^N\tk,Y\NooI\tk,Z^{N,I}\tk\rp\rc\rv^2\\
&+\lp1+\frac{ h}{\e}\rp \E\lv P_k\lc\DeltaB{k}\gTNIkk\rc\rv^2\\
\leq 
&\lp1+\frac{ h}{\e}\rp \E\lv \Ek Y^{N,I,I}\tkk \rv^2 
+2(h^2 +2\e h)\lv  f(0,0,0)\rv^2\\
&+2L_f(h^2 +2\e h)\lp\E\lv  X^N\tk\rv^2+\E\lv  Y\NooI\tk\rv^2+\E\lv Z^{N,I}\tk\rv^2 \rp\\
&+\lp1+\frac{ h}{\e}\rp \E\lv\Ek\lc\DeltaB{k}\gTNIkk\rc\rv^2.
\end{align*}
Using the upper estimate \eqref{bound-ZNI} in Lemma \ref{lem-2}, we obtain
\begin{align*}
[ 1-& 2L_f( h^2+2\e h)]\E\lv Y^{N,\infty,I}\tk\rv^2\\
\leq &
 \lp1+\frac{ h}{\e}-2L_f\lp h+2\e\rp  \rp \E\lv \Ek Y^{N,I,I}\tkk  \rv^2
 +2\lp h^2+\e h\rp \lp|f(0,0,0)|^2+L_f\E\lv X^N\tk\rv^2\rp\\
&+2L_f\lp h+2\e\rp  \E\lv  Y^{N,I,I}\tkk \rv^2 
+2L_f\lp h+2\e\rp \E\lv \DeltaB{k} \gTNIkk\rv^2 \\
&+\lp1+\frac{ h}{\e}-2L_f\lp h+2\e\rp \rp\E\lv \E_k\lp\DeltaB{k} \gTNIkk\rp\rv^2 .
\end{align*}
Choose $\e$ such that $4L_f\e=1$. Then $\lp 1+\frac{h}{\e}\rp-2L_f\lp h+2\e\rp=2L_fh$
and $2L_f( h+2\e )=2L_f h+ 1$. 
Using Theorem \ref{lem-euler-X} we deduce the exitence of
$C>0$ such that,
\begin{align*}
 \lc 1- 2L_f( h^2+2\e h)\rc&\E\lv Y^{N,\infty,I}\tk\rv^2\\
\leq & Ch+\lp 1+4L_fh\rp\lc  \E\lv  Y^{N,I,I}\tkk \rv^2 
+ \E\lv \DeltaB{k} \gTNIkk\rv^2 \rc.
\end{align*}
Then for $h^*\in(0,1]$ small enough (ie $(2L_f+1)h^*<1$),
using Lemma \ref{lem-1-Ch}, we deduce that for
 $\Gamma:=\frac{2L_f+1}{1-(2L_f+1)h^*}$
and $h\in(0,h^*)$, we have $(1-(2L_f+1)h)^{-1}\leq 1+\Gamma h$.
Thus using the independence of $\DeltaB{k}$ and $\F\tkk$,
the growth condition \eqref{cg} and Lemma \ref{lem-euler-X}, we
deduce the existence of a constant $C>0$, 
such that for $h\in(0,h^*)$,
\begin{align*}
\E\lv Y^{N,\infty,I}\tk\rv^2
\leq & Ch+\lp 1+Ch\rp  \E\lv  Y^{N,I,I}\tkk \rv^2 
+C h \E\lv  \gTNIkk\rv^2\\
\leq & Ch+\lp 1+Ch\rp  \E\lv  Y^{N,I,I}\tkk \rv^2 .
\end{align*}
This concludes the proof of \eqref{bound-Yinf}.
\endproof
Let $\eta_k^{N,I}:=\E\lv Y\NII\tk-Y^N\tk\rv^2$ for $k=0,\dots,N$; the following lemma gives an
upper bound of the $L^2$-norm of $Y\NinftyI\tk-P\k Y^N\tk$ in terms of 
$\eta\NI\kk$.
\begin{lemma}\label{lem-25}
For small $h$ and for $k=0,\dots,N-1$ we have:
\begin{align*}
\E\lv\Yinf{k}-P_k\YNt{k}\rv^2
\leq(1+K h) \eta^{N,I}_{k+1}
+  Kh\lc\E\lv R_k\YNt{k}\rv^2 
+  \E\lv R_kZ^N\tk\rv^2\rc .
\end{align*}
\end{lemma}
\proof
The argument, which is similar to that in the proof of 
Lemmas \ref{lem-2} and \ref{lem-3-8} is more briefly sketched.
Applying the operator $P_k$ to both sides of equation \eqref{equationYN} and using 
 \eqref{lem-PE=P}, we obtain
\begin{align*}
P_kY^N\tk
=&P_k Y^N\tkk+hP_k\lc f\XYZNt{k}\rc +P_k\lc\DeltaB{k}\gTNIkk\rc.
\end{align*}
Hence Lemma \ref{lem-3-8} implies that
\begin{align*}
\Yinf{k}-P_k\YNt{k}=&P_k\lc\YNIt{I}{k+1}-\YNt{k+1}\rc +hP_k\lc f\XYZinf{k}-f\XYZNt{k}\rc\\
&+P_k\lp \DeltaB{k}\lc\gTNIkk-\gXYNkk\rc\rp.
\end{align*}
Lemma \ref{lem-1} proves the orthogonality of the first and third term of the above decomposition.
Squaring this equation, using Young's inequality and \eqref{lem-PE=P}, 
the $L^2$-contraction property of $P\k$ and the Lipschitz property of $g$ given in \eqref{cg},
computations similar to that made in the proof of Lemma \ref{lem-3-8} yield
\begin{align}
\E&\lv\Yinf{k}-P_k\YNt{k}\rv^2= \lp1+\frac{h}{\e}\rp\E\lv P_k\lc\YNIt{I}{k+1}-\YNt{k+1}\rc \rv^2\nonumber\\
&+h^2\lp 1+2\frac{\e}{h}\rp\E\lv P_k\lc f\XYZinf{k}-f\XYZNt{k}\rc\rv^2
\nonumber\\
&+\lp1+\frac{h}{\e}\rp\E\lv P_k\lp \DeltaB{k}\lc\gTNIkk-\gXYNkk \rc\rp\rv^2
\nonumber\\
\leq &
 \lp1+\frac{h}{\e}\rp\E\lv \Ek\lc\YNIt{I}{k+1}-\YNt{k+1}\rc \rv^2
+L_f\lp h+2\e\rp h\lp\E\lv  \Yinf{k}-\YNt{k}\rv^2+\E\lv Z^{N,I}\tk-Z^N\tk\rv^2\rp\nonumber\\
&+\lp1+\frac{h}{\e}\rp\E\lv \Ek\lp \DeltaB{k}\lc\gTNIkk-\gXYNkk\rc \rp\rv^2.\label{ineq-24}
\end{align}
By construction $\Yinf{k}\in\P\k$. Hence
\begin{equation}\label{eq-22}
\E\lv  \Yinf{k}-\YNt{k}\rv^2= \E\lv  \Yinf{k}-P_k\YNt{k}\rv^2+\E\lv R_k\YNt{k}\rv^2.
\end{equation}
Using Lemma \ref{lem-2} 
we deduce that for any $\epsilon>0$
\begin{align*}
&\lp 1-L_f\lp h^2+2\e h\rp\rp \E\lv\Yinf{k}-P_k\YNt{k}\rv^2\\
\leq& L_f\lp h+2\e\rp \eta^{N,I}_{k+1}
+hL_f\lp h+2\e\rp\lc \E\lv R_k\YNt{k}\rv^2 
+\E\lv R_kZ^N\tk\rv^2\rc\\
&+\lp \lp1+\frac{h}{\e}\rp- L_f\lp h+2\e\rp \rp \E\lv\Ek\lp \YNIt{I}{k+1}-Y^N\tkk\rp\rv^2 \\
&+ L_f\lp h+2\e\rp\E\lv\DeltaB{k}\lc\gTNIkk-\gXYNkk\rc \rv^2\\
&+\lp \lp1+\frac{h}{\e}\rp-L_f\lp h+2\e\rp  \rp 
\E\lv\Ek\lp \DeltaB{k}\lc\gTNIkk
-\gXYNkk\rc \rp\rv^2.
\end{align*}
Let $\e>0$ satisfy $2L_f\e=1$; then 
$\lp1+\frac{h}{\e}\rp-L_f\lp h+2\e\rp =L_fh$ 
and $L_f\lp h+2\e \rp=L_fh+1$.
Thus, since $\E\k$ contracts the $L^2$-norm, 
we deduce
\begin{align*}
&\lp 1-L_f\lp h^2+2\e h\rp\rp \E\lv\Yinf{k}-P_k\YNt{k}\rv^2\\
\leq& \lp 1+2L_fh\rp \eta^{N,I}_{k+1}
+h\lp 1+L_fh\rp \lc\E\lv R_k\YNt{k}\rv^2 
+\E\lv R_kZ^N\tk\rv^2\rc\\
&+ \lp 1+2L_fh\rp\E\lv\DeltaB{k}\lc\gTNIkk-\gXYNkk\rc \rv^2.
\end{align*}
Let $h^*\in(0,\frac{1}{L_f+1}$) and set
$\Gamma=\frac{L_f+1}{1-(L_f+1)h^*}$. Lemma \ref{lem-1-Ch} shows that for $h\in(0,h^*)$
we have $\lp 1-L_f\lp h^2+2\e h\rp\rp^{-1}\leq 1+\Gamma h$.
The previous inequality, the independence of 
$\DeltaB{k}$ and $\F\tkk$ and the Lipschitz property \eqref{cg}
imply that for some constant $K$ which can change for one line to the next
\begin{align*}
\E&\lv\Yinf{k}-P_k\YNt{k}\rv^2
\leq(1+K h) \eta^{N,I}_{k+1}
+  Kh\lc\E\lv R_k\YNt{k}\rv^2 
+  \E\lv R_kZ^N\tk\rv^2\rc\\
&+K h\E\lv\gTNIkk-\gXYNkk \rv^2\\
\leq&(1+K h) \eta^{N,I}_{k+1}
+  Kh\lc\E\lv R_k\YNt{k}\rv^2 
+\E\lv R_kZ^N\tk\rv^2\rc.
\end{align*}
This concludes the proof of Lemma \ref{lem-25}
\endproof
The following Lemma provides $L^2$-bounds of $Y\NII\tk$,
$Y\NinftyI\tk$ and $Z\NI\tk$ independent of $N$ and $I$.
\begin{lemma}\label{bound-YYY}
There exists a constant $K$ such that  for
large $N$ and for every $I\geq1$, 
$$
\max_{0\leq k\leq N}\E\lv Y^{N,I,I}\tk\rv^2+
\max_{0\leq k\leq N-1} \E\lv Y^{N,\infty,I}\tk\rv^2
+\max_{0\leq k \leq N} h\E\lv Z^{N,I}\tk\rv^2\leq K.
$$
\end{lemma}
\proof
Using inequality \eqref{error-Yinf-YNiI} and Young's inequality, we have the following bound, 
for $i=1,\dots,I$, $h<1$ and some constant $K$ depending on $L_f$:
\begin{align}
\E\lv Y^{N,i,I}\tk\rv^2\leq & \lp1+\frac{1}{h}\rp\E\lv \Yinf{k}-Y^{N,i,I}\tk\rv^2+(1+h)\E\lv \Yinf{k}\rv^2\nonumber\\
\leq & \lp1+\frac{1}{h}\rp L_f^ih^{2i}\E\lv \Yinf{k}\rv^2+(1+h)\E\lv \Yinf{k}\rv^2
\leq  (1+Kh)\E\lv \Yinf{k}\rv^2.\label{eqq1}
\end{align}
Choosing $i=I$ and using \eqref{bound-Yinf}
we deduce that for some constant $K$ which can change from line to line,
$\E\lv Y^{N,I,I}\tk\rv^2\leq  Kh+(1+Kh)\E\lv Y^{N,I,I}\tkk \rv^2$.
Hence Lemma \ref{lem-Gronwall} yields $\max_k\E\lv Y^{N,I,I}\tk\rv^2\leq  K$.
Plugging this relation into inequality \eqref{bound-Yinf} proves that
$$
\max_k\E\lv Y\NII\tk\rvd+\max_k \E\lv Y^{N,\infty,I}\tk\rv^2\leq  K<\infty.
$$
Using \eqref{bound-ZNI} and the independence of $\DeltaB{k}$ and 
$\F\tkk$, we deduce
\begin{align*}
h\E\lv Z^{N,I}\tk\rv^2
\leq& \E\lv Y\NII\tkk\rv^2+\E\lv \DeltaB{k}g\lp X^N\tkk,Y\NII\tkk\rp\rv^2\\
\leq  & \E\lv Y\NII\tkk\rv^2+h\E\lv g\lp X^N\tkk,Y\NII\tkk\rp\rv^2
\end{align*}
Finally, the Lipschitz property \eqref{cg}
yields
\begin{align*}
h\E\lv Z^{N,I}\tk\rv^2
\leq  & \E\lv Y\NII\tkk\rv^2+2h\lv g\lp 0,0\rp\rv^2 +2hL_g \lp\E\lv  X^N\tkk\rv^2+\E\lv Y\NII\tkk\rv^2 \rp\\
\leq  &\lp 1+2hL_g\rp \E\lv Y\NII\tkk\rv^2+2h\lv g\lp 0,0\rp\rv^2 +2hL_g \E\lv  X^N\tkk\rv^2.
\end{align*}
Theorem \ref{lem-euler-X} and the $L^2$-upper estimates of $Y\NII\tkk$ conclude the proof.
\endproof
The following lemma provides a backward recursive upper estimate of $\eta\NI_.$
Recall that $\eta\NI\k=\E\lv Y\NII\tk-Y^N\tk\rvd$
\begin{lemma}\label{lem-20}
For $0\leq k < N$, we have:
$$
\eta^{N,I}_k\leq  (1+K h) \eta^{N,I}_{k+1}+Ch^{2I-1}
+  K\E\lv R_k\YNt{k}\rv^2 
+ Kh \E\lv R_kZ^N\tk\rv^2.
$$
\end{lemma}
\proof
For $k=N$, $Y^N_{t_N}=\Phi\lp X^N_{t_N}\rp$ 
and $Y\NII_{t_N}=P_N\Phi\lp X^N_{t_N}\rp$ so that 
$\eta^{N,I}_N=\E\lv \Phi\lp X^N_{t_N}\rp-P_N \Phi\lp X^N_{t_N}\rp\rvd$. 
Let $k\in\{0,\dots, N-1\}$; using inequality \eqref{error-Yinf-YNiI} and Young's inequality,
we obtain
\begin{align*}
\eta^{N,I}_k=&\E\lv Y\NII\tk-Y^N\tk\rv^2\\
\leq &\lp 1+\frac{1}{h}\rp \E \lv Y\NII\tk-Y^{N,\infty,I}\tk\rv^2+(1+h)\E\lv Y^{N,\infty,I}\tk - Y^N\tk\rv^2\\
\leq & \lp 1+\frac{1}{h}\rp L_f^Ih^{2I}\E \lv Y^{N,\infty,I}\tk\rv^2 
+(1+h)\E\lv Y^{N,\infty,I}\tk - P_k Y^N\tk\rv^2+(1+h)\E\lv R_k Y^N\tk\rv^2.
\end{align*}
Finally, Lemmas \ref{bound-YYY} and \ref{lem-25} imply that for
some constant $K$ we have for every $N$ any $k=1,\dots,N$:
\begin{align}
\eta^{N,I}_k\leq & Kh^{2I-1} +(1+h)\E\lv Y^{N,\infty,I}\tk - P_k Y^N\tk\rv^2+(1+h)\E\lv R_k Y^N\tk\rv^2\label{ineq-21}\\
\leq & (1+K h) \eta^{N,I}_{k+1}+Kh^{2I-1}
+  K\E\lv R_k\YNt{k}\rv^2 
+ Kh \E\lv R_kZ^N\tk\rv^2;\nonumber
\end{align}
this concludes the proof. 
\endproof
Gronwall's Lemma \ref{lem-Gronwall} and Lemma \ref{lem-20} prove the existence of $C$ such that 
for $h$ small enough
\begin{align}
\max_{0\leq k\leq N}\E\lv\YNIt{I}{k}-\YNt{k}\rv^2
\leq & C h^{2I-2}+C\sum_{k=0}^{N-1}\E\lv R_k Y^N\tk\rv^2
+Ch\sum_{k=0}^{N-1}\E\lv R_k Z^{N,I}\tk\rv^2\nonumber\\
& + C \E\lv \Phi\lp X^N_{t_N}\rp - P_N\Phi\lp X^N_{t_N}\rp\rvd\label{4.15Bis}
\end{align}
which is part of Theorem \ref{th-2}.
Let $\zeta^N:=h\sum_{k=0}^{N-1}\E\lv Z^{N,I}\tk-Z^N\tk\rv^2$. 
In order to conclude the proof Theorem \ref{th-2}, 
we need to upper estimate $\zeta^N$, 
which is done in the next lemma.
\newcommand{\zetaN}{\zeta^N}
\newcommand{\sumkN}{\sum_{k=0}^{N-1}}
\newcommand{\ImIM}{^{I-1,I,M}}
\begin{lemma}\label{lem-der}
There exits a constant $C$ such that for $h$ small enough
and every $I\geq1$
\begin{align*}
\zetaN\leq C h^{2I-2}+Ch\sumkN \E\lv R_k Z^N_k\rv^2
+C\sumkN \E\lv R_k Y^N_k\rv^2+C\max_{0\leq k\leq N-1}\eta^{N,I}_k.
\end{align*}
\end{lemma}
\proof
Multiply inequality \eqref{bound-ZNI-ZN} by $h$,  use the independence of $\DeltaB{k}$
and $\F\tkk$ and the Lipschitz property \eqref{cg}; this yields 
\begin{equation}\label{ineq-sur-zeta}
\zetaN\leq h\sumkN\E\lv R_kZ^N\tk\rv^2
+\sumkN\lp \lp 1+L_gh\rp\E\lv \YNIt{I}{k+1}-Y^N\tkk\rv^2 -\E\lv\Ek\lp \YNIt{I}{k+1}-Y^N\tkk\rp\rv^2 \rp.
\end{equation}
Multiply inequality\eqref{ineq-24} by $(1+L_gh)(1+h)$,
use the independence of $\DeltaB{k}$ and $\F\tkk$
and the Lipschitz property \eqref{cg}; this yields for $\e>0$:
\begin{align}
(1+L_gh)&(1+h)\E\lv\Yinf{k}-P_k\YNt{k}\rv^2
\nonumber\\
\leq &
 \lp1+\frac{h}{\e}\rp(1+L_gh)(1+h)\E\lv \Ek\lc\YNIt{I}{k+1}-\YNt{k+1}\rc \rv^2\nonumber\\
&+L_f\lp h+2\e\rp h(1+L_gh)(1+h)\lp\E\lv  \Yinf{k}-\YNt{k}\rv^2+\E\lv Z^{N,I}\tk-Z^N\tk\rv^2\rp\nonumber\\
&+\lp1+\frac{h}{\e}\rp(1+L_gh)(1+h)L_gh\E\lv \YNIt{I}{k+1}-\YNt{k+1} \rv^2.\label{eq-page-31}
\end{align}
Multiply inequality \eqref{ineq-21} by $\lp 1+L_gh\rp$ and use \eqref{eq-page-31}; 
this yields for some constants $K$, $C$, $\bar C$ and $h\in(0,1]$, $\e>0$:
\begin{align*}
\Delta\kk:=& \lp 1+L_gh\rp\E\lv \YNIt{I}{k+1}-Y^N\tkk\rv^2 -\E\lv\Ek\lp \YNIt{I}{k+1}-Y^N\tkk\rp\rv^2 \\
\leq & 
  K  h^{2I-1} +K\E\lv R_k Y^N\tk\rv^2
+\lp \lp1+\frac{h}{\e}\rp(1+L_gh)(1+h)-1\rp\E\lv \Ek\lc\YNIt{I}{k+1}-\YNt{k+1}\rc \rv^2\\
&+C\lp h+2\e\rp h\lp\E\lv  \Yinf{k}-\YNt{k}\rv^2+\E\lv Z^{N,I}\tk-Z^N\tk\rv^2\rp\\
&+\lp1+\frac{h}{\e}\rp Ch\E\lv \YNIt{I}{k+1}-\YNt{k+1} \rv^2.
\end{align*}
Now we choose $\e$ such that $2C\e=\frac{1}{4}$; 
then we have for some constant $K$ and $h\in(0,1]$:
\begin{align*}
\Delta\kk
\leq & 
  K  h^{2I-1} +K\E\lv R_k Y^N\tk\rv^2+Kh\E\lv \YNIt{I}{k+1}-\YNt{k+1} \rv^2\\
&+\lp Ch+\frac{1}{4}\rp h\lp\E\lv  \Yinf{k}-\YNt{k}\rv^2+\E\lv Z^{N,I}\tk-Z^N\tk\rv^2\rp.
\end{align*}
Thus, for $h$ small enough (so that $Ch\leq \frac{1}{4}$), summing over $k$
we obtain
\begin{align*}
\sumkN\Big( \lp 1+L_gh\rp&\E\lv \YNIt{I}{k+1}-Y^N\tkk\rv^2 -\E\lv\Ek\lp \YNIt{I}{k+1}-Y^N\tkk\rp\rv^2\Big) \\
\leq & 
  K  h^{2I-2} +K\sumkN\E\lv R_k Y^N\tk\rv^2+K\max_k\eta^{N,I}\k\\
&+\undeux h\sumkN\lp\E\lv  \Yinf{k}-\YNt{k}\rv^2+\E\lv Z^{N,I}\tk-Z^N\tk\rv^2\rp.
\end{align*}
Plugging this inequality in \eqref{ineq-sur-zeta} yields
\begin{align*}
\undeux\zetaN\leq& 
 K  h^{2I-2} +h\sumkN\E\lv R_kZ^N\tk\rv^2+K\sumkN\E\lv R_k Y^N\tk\rv^2+K\max_k\eta^{N,I}\tk\\
&+\undeux h\sumkN\E\lv  \Yinf{k}-\YNt{k}\rv^2.
\end{align*}
Using \eqref{eq-22} and Lemma \ref{lem-25}, we obtain for some constant $K$ and every $h\in(0,1]$
\begin{align*}
h\sumkN& \E\lv  \Yinf{k}-\YNt{k}\rv^2
\leq
(1+K h) h\sumkN\eta^{N,I}_{k+1}+  Kh^2\sumkN\lc\E\lv R_k\YNt{k}\rv^2 
+ \E\lv R_kZ^N\tk\rv^2\rc\\
&+h\sumkN\E\lv R_k\YNt{k}\rv^2
\leq
K \max_k\eta^{N,I}_{k}+  K\sumkN\E\lv R_k\YNt{k}\rv^2 
+ Kh\sumkN \E\lv R_kZ^N\tk\rv^2.
\end{align*}
This concludes the proof of Lemma \ref{lem-der}.
\endproof
Theorem \ref{th-2} is a straightforward consequence of inequality \eqref{4.15Bis} and Lemma \ref{lem-der}.
\section{Approximation step 3}
In this section we will use regression approximations
and introduce some minimization problem for a 
$M$-sample of $(B,W)$ denoted by $(B^m,W^m,m=1,\dots,M)$.
This provides a Monte Carlo approximation of $Y^N,I,I$ and
$Z^N,I$ on the time grid.
\subsection{Some more notations for the projection}
\it{We at first introduce some notations}
\begin{itemize}
\item[\bf{(N5)}] \it{For fixed $k=1,\dots,N$ and $m=1,\dots,M$, let $p_k^m$ denote the  orthonormal 
family of $L^2\lp \Omega\rp$
similar to $p_k$ in (N4) replacing $X^N$ by $X^{N,m}$
and $B$ by $B^m$.}
\item[\bf{(N6)}] 
\it{For a real $n\times n$ symmetric matrix $A$, $\| A\|$ is the maximum of the absolute value of its eigenvalues 
and $\|A\|_F=\lp\sum_{i,j}A_{i,j}^2\rp^{\undeux}$ its Frobenius norm. 
If 
$A:\R^n\rightarrow\R^n$ also denotes the linear operator whose matrix in the canonical basis is $A$, then 
$\|A\|$ is the operator-norm of $A$ when $\R^n$ is endowed with the Euclidian norm.
Note that $\|A\|\leq \|A\|_F$ follows from Schwarz's inequality.}
\item[\bf{(N7)}] 
\it{For $k=0,\dots,N-1$ and $m=1,\dots,M$ 
let  $ v^m_k$ and $ v_k$ be  column vectors
whose entries are the components in the canonical base of the vectors
\begin{align}
\lp p\mk,p\mk\frac{\DeltaWm{k+1}}{\sqrt h}\rp,\;and\;
\lp p_k,p_k\frac{\DeltaW{k+1}}{\sqrt h}\rp\label{def-vmk}
\end{align}
respectively.
Note that $\E v\k v\k^*=Id$,
since the entries of $p\k$ are an orthonormal family of 
$L^2\lp\F\k\rp$ and $\frac{\Delta W\kk}{h}$
is a normed vector in $L^2$ independent of $p_k$.}
\item[\bf{(N8)}] 
\it{For $k=0,\dots,N-1$ let 
$V\Mk,P^M_k$ be symmetric matrices defined by 
\begin{align}
V\Mk:=\unM\summ v\mk \lc v^m_k\rc^*,
P^M_k:=\unM\summ p\mk (p\mk)^*.
\label{def-VMk}
\end{align}
}
\item[\bf{(N9)}] 
\it{We denote by $\mathcal N$ the 
$\sigma$-algebra of measurable sets $A$ with $\mathbb P(A)=0$ and set:
\begin{align*}
\F^{W,m}_{t}=&\sigma \lp W^m_s ; 0\leq s\leq t\rp\vee\mathcal N, \quad
\F^{B,m}_{t,t'}= \sigma \lp B^m_s-B^m_{t'}; t\leq s\leq t'\rp\vee\mathcal N,\\
\F^{W,M}_{t}=&\F^{W}_{t}\vee\bigvee_{m=1}^M\F^{W,m}_{t},\quad
\F^{B,M}_{t,T}= \F^{B}_{t,T}\vee\bigvee_{m=1}^M\F^{B,m}_{t,T},\quad
\F_t=\F^W_t\vee\F^B_{t,T}.
\end{align*}
Note that $\lp\F_t\rp_t$ and $\lp\F^B_{t,T}\rp_t$ are not filtrations.}
\item[\bf{(N10)}] 
\it{In the sequel we will need to localize some processes using the following events 
\begin{align}
\mathfrak A_j:=&\left\lbrace \| V^M_j-Id\|\leq h ,\|P^M_j-Id \|\leq h\right\rbrace
\in \F^{W,M}_{t_{j+1}}\vee\F^{B,M}_{t_j,T},\label{def-mathfrak-Aj}\\
\AMk:=&\bigcap_{j=k}^{N-1} \mathfrak A_j\in \F^{W,M}_{t_{N}}\vee\F^{B,M}_{t_k,T}.\label{def-AMk}
\end{align}
}
\item[\bf{(N11)}]
\it{For $x=(x_1,\dots,x_m)\in\R^M$,
we denote $\lv x\rvd_M:=\unM\summ\lv x_m\rvd$.}
\end{itemize}
\subsection{Another look at the previous results }
We introduce the following random variables
\begin{align*}
\zeta^N_k:=\rho^N_k:=\lp \lv p\k\rv\sqrt{C_0}\rp\vee 1,
\end{align*}
where $C_0$ is  constant in the Lemma \ref{bound-YYY}.
Since $\YNIt{i}{k}$ and $Z^{N,I}\tk$ are in $\P_k$ (see \eqref{eq-ZNI} and \eqref{eq-YNiI}),
we can rewrite these random variables as follows:
\begin{align}
Y\NiI\tk=\alpha^{i,I}_k .p_k=\lp\alpha^{i,I}_k\rp^* p_k,\quad 
Z^{N,I}\tk=\beta^I_k. p_k=\lp\beta^I_k\rp^*pk, \label{eq-alpha-beta}
\end{align}
where $\alpha^{i,I}_k$ (resp. $\beta^I\k$) is the vector of the coefficient 
in the basis $p\k$ of the random variable $Y\NiI\tk$ (resp. $Z\NI\tk$),
identified with the column matrix of the coefficients in the canonical basis.
\begin{remark}Note that the vectors $\alpha^{i,I}_k$ and $\beta^I_k$ are deterministic.
\end{remark}
The following Proposition gives a priori estimates of 
$Y^{N,i,I}\tk$ and $Z^{N,I}\tk$.
\begin{proposition}\label{a-priori-estimate}
For $i\in\lbrace 1,\dots,I\rbrace\cup\{\infty\}$
and for $k=0,\dots,N$, 
we have $\lv Y^{N,i,I}\tk\rv\leq \rho^N_k,\quad \sqrt h \lv Z\NI\tk\rv\leq\zeta\Nk.$
Moreover, for every $I$ and $i=0,\dots,I$:
\begin{equation}\label{lem-alpha-beta}
\lv\alpha\iI\k\rvd\leq\E\lv\rho\Nk\rvd,\quad \lv\beta^I\k\rvd\leq\unh\E\lv\zeta\Nk\rvd.
\end{equation}
\end{proposition}
\proof
Let $i\in\{1,\dots,I\}\cup\lbrace \infty\rbrace$ 
and $k=0,\dots,N$.
Squaring $Y\NiI\tk$,
taking expectation and using the previous remark, we obtain
\begin{align*}
\E\lv Y\NiI\tk\rvd=\lp\alpha^{i,I}_k\rp^*\E\lp p_k  p\k^*\rp\alpha^{i,I}_k
\geq\lp\alpha^{i,I}_k\rp^*\alpha^{i,I}_k
=\lv\alpha^{i,I}_k\rvd
\end{align*}
Using Lemma \ref{bound-YYY}, we deduce that
$\lv\alpha^{i,I}_k\rvd\leq C_0.$
The Cauchy-Schwarz inequality implies 
\begin{align*}
\lv Y\NiI\tk\rv\leq\lv\alpha^{i,I}_k\rv\lv p_k\rv
\leq \lv p\k\rv\sqrt C_0
\leq \lp\lv p\k\rv\sqrt C_0\rp\vee 1.
\end{align*}
A similar computation based on Lemma \ref{bound-YYY} proves that $\sqrt h \lv Z\NI\tk\rv\leq \zeta\Nk.$
The upper estimates of $\lv \alpha^{i,I}\k\rvd$ and $\lv\beta^I\k\rvd$ are straightforward consequences
of the previous ones.
\endproof
We now prove that $\lp\alpha^{i,I}\k,\beta^I\k\rp$ solves a minimization problem.
\begin{proposition} 
The vector $\lp\alpha^{i,I}_k,\beta^I_k\rp$ solves
the following minimization problem:
for $k=0,\dots,N-1$ and
for every $i=1,\dots,I$, we have:
\begin{align}
\lp\ay{i}{k},\az{I}{k}\rp=\arg\min_{\lp\alpha,\beta\rp}\E\Big|&Y^{N,I,I}_{k+1}-\alpha.\py{k}
+hf\lp X^N\tk,\ay{i-1}{k}.p_k,Z^{N,I}\tk\rp\nonumber\\
&+\DeltaB{k} g\lp \XNt{k+1},Y\NII\tkk\rp- \beta.\pz{k}\Delta W_{k+1}\Big|^2.\label{eq-9}
\end{align}
\end{proposition}
\proof 
Let $(Y,Z)\in \P_k\times \P_k$; then since $\P\k\subset L^2\lp\mathcal F\tk\rp$ and $\Delta W\kk$ is independent of 
$\mathcal F\tk$, we have
\begin{align*}
 \E&\lv Y\NII\tkk-Y+hf\lp\XNt{k},Y^{N,i-1,I}\tk,Z\tk^{N,I}\rp+\DeltaB{k} g\lp \XNt{k+1},Y\NII\tkk\rp-Z\DeltaW{k+1}\rv^2\\
=&\E\lv Y\NII\tkk-Y+hf\lp \XNt{k},Y^{N,i-1,I}\tk,Z\tk^{N,I}\rp +\DeltaB{k} g\lp \XNt{k+1},Y\NII\tkk\rp\rv^2\\
&+h\E\lv Z-\frac{1}{h}\lp Y\NII\tkk+\DeltaB{k} g\lp\XNt{k+1},Y\NII\tkk\rp \rp\DeltaW{k+1}\rv^2\\
&-\frac{1}{h}\E\lv \lp Y\NII\tkk+\DeltaB{k} g\lp\XNt{k+1},Y\NII\tkk\rp\rp\DeltaW{k+1}\rv^2.
\end{align*}
The minimun on pairs of elements of $\P\k$ is given by 
the orthogonal projections, that is by the random variables
$Y=Y\NiI\tk$ and $Z=Z\NI\tk$ defined by
\eqref{eq-YNiI} and \eqref{eq-ZNI} respectively. This concludes the proof
using the notations introduced in \eqref{eq-alpha-beta}.
\endproof
For $i\in\lb1,\dots,I\rb\cup \lb \infty\rb$, we define 
$\theta\iI\k:=\lp \alpha\iI\k,\sqrt h\beta^I\k\rp$.
The following lemma gives some properties on $\theta\iI\k$.
\begin{lemma}\label{lem-33}
For all $i\in\lbrace1,\dots,I\rbrace\cup\lbrace\infty\rbrace$, we have 
for $k=0,\dots,N$
(resp. for $k=0,\dots,N-1$)
\begin{align*}
\lv \theta\iI\k\rv^2\leq  \E \lv\rho \Nk\rv ^2+\E \lv \zeta\Nk\rv^2,\quad
resp.\quad
\lv \theta\inftyI\k-\theta\iI\k\rv^2\leq L_f^ih^{2i}\E\lv\rho\Nk\rv^2.
\end{align*}
Furthermore, we have the following explicit expression of $\theta\inftyI_k$ for $v\k$ defined by \eqref{def-vmk}:
\begin{align}\label{eq-34}
\theta\inftyI_k=\E\lc v_k \lp\alpha\II\kk .p\kk+hf\lp X\Nk,\alpha\inftyI\k. p\k,\beta^I\k. p\k\rp 
+\DeltaB{k} g\lp X^N\tkk,\alpha\II\kk. p\kk\rp\rp\rc.
\end{align}
\end{lemma}
\proof
Proposition \ref{a-priori-estimate} implies that
$
\lv \theta\iI\k\rv^2=\lv \alpha\iI\k\rv^2+h\lv \beta^I\k\rv^2
\leq  \E \lv\rho \Nk\rv ^2+\E \lv \zeta\Nk\rv^2.
$\\
Using inequality \eqref{error-Yinf-YNiI} and  Proposition \ref{a-priori-estimate}, since $\E\lv p\k\rvd=1$ we obtain 
\begin{align*}
\lv \theta\inftyI\k-\theta\iI\k\rv^2
=\E \lv Y\NinftyI\tk-Y\NiI\tk\rv^2
\leq L_f^ih^{2i}\E\lv Y\NinftyI\tk\rv^2
\leq L_f^ih^{2i}\E\lv \rho\Nk\rv^2.
\end{align*}

Using equation \eqref{eq-Yinf} and the fact that the components of $p\k$ are  an orthonormal family of $L^2$,
we have
\begin{align*}
\alpha\inftyI\k=&\E\lc p\k Y\NinftyI\k\rc\\
=&\E \lp p\k  P_{k}\lc Y^{N,I,I}\tkk+hf\lp X^N\tk,Y\NooI\tk,Z^{N,I}\tk\rp +\DeltaB{k}\gTNIkk\rc\rp\\
=&\E \lc p\k\lp\alpha\II\kk .p\kk+hf\lp X\Nk,\alpha\inftyI\k .p\k,\beta^I\k .p\k\rp 
+\DeltaB{k} g\lp X^N\tkk,\alpha\II\kk .p\kk\rp\rp\rc.
\end{align*}
A similar computation based on equation \eqref{eq-ZNI} 
and on the independence of $\F\tk$ and $\Delta W\kk$ yields
\begin{align*}
\sqrt h\beta^I\k=&\E \lc\sqrt h p\k Z\NI\tk \rc\\
=&\E\lc  \frac{1}{\sqrt h} p\k P\k\lp    Y^{N,I,I}\tkk\Delta W_{k+1}+ \DeltaB{k}\gTNIkk\DeltaW{k+1}\rp\rc\\
=&\E\lc p\k\lp\alpha\II\kk .p\kk\frac{\Delta W\kk}{\sqrt h}
+\DeltaB{k} g\lp X^N\tkk,\alpha\II\kk .p\kk\rp\frac{\Delta W\kk}{\sqrt h}\rp\rc\\
=&\E\lc p\k\frac{\Delta W\kk}{\sqrt h}\lp\alpha\II\kk .p\kk
+h f\lp X\Nk,\alpha\inftyI\k .p\k,\beta^I\k .p\k\rp 
+\DeltaB{k} g\lp X^N\tkk,\alpha\II\kk .p\kk\rp \rp\rc.
\end{align*}
Finally, we recall by \eqref{def-vmk} that $v\k:=\lp p\k,p\k\frac{\Delta W\kk}{\sqrt h}\rp$; this concludes the proof.
\endproof
\subsection{The numerical scheme}
Let $\xi:\R\rightarrow\R$ be a $C^2_b$ function, such that 
$\xi(x)=x$ for $|x|\leq 3/2, |\xi|_{\infty}\leq 2$ and $|\xi'|_{\infty}\leq 1$.
We define the random truncation functions
\begin{align}\label{rhozeta}
\widehat\rho_k^N(x):=\rho^N_k\xi\lp\frac{x}{\rho^N_k}\rp,  \quad
\widehat\zeta_k^N(x):=\zeta^N_k\xi\lp\frac{x}{\zeta^N_k}\rp  .
\end{align}
The following lemma states some properties of these functions.
\begin{lemma}\label{lem-inv}
Let $\widehat\rho\k^N$ and $\widehat\zeta\Nk$ be defined by \eqref{rhozeta}, then
  \begin{enumerate}
  \item $\wr$ (resp. $\wz$) leaves $Y\NII\tk$ (resp. $\sqrt hZ\NI\tk$) invariant, that is:
\begin{align*}
  \widehat\rho^N_k\lp \alpha^{I,I}_k.p_k\rp=\alpha^{I,I}_k.p_k,\quad
  \widehat\zeta\Nk\lp\sqrt h\beta_k^{I}. p_k\rp=\sqrt h\beta_k^{I}. p_k.
\end{align*}
\item $\wr,\wz$ are 1-Lipschitz and 
$\lv\widehat \rho\Nk(x)\rv\leq \lv x\rv$ for every $x\in\R$.
\item  $\wr$ (resp. $\wz$) is bounded by $2\lv \rho^N_k\rv$ (resp. by  $2\lv \zeta^N_k\rv$).
 \end{enumerate}
\end{lemma}
\proof 
In part (1)-(3) we only give the proof for $\widehat \rho\Nk$, since that for $\widehat\zeta\Nk$ is similar.\\
1.
By Proposition \ref{a-priori-estimate}, $\lv\frac{ \alpha\II\k. p\k}{\rho\Nk}\rv\leq 1$. 
Hence, $\xi\lp\frac{ \alpha\II\k .p\k}{\rho\Nk}\rp=\frac{ \alpha\II\k .p\k}{\rho\Nk}$.\\
2. Let $y,y'\in \R$; since $\lv \xi'\rv_{\infty}\leq1$,
\begin{align*}
\lv\wr(y)-\wr(y')\rv
=\lv\rho^N_k\rv\lv\xi\lp\frac{y}{\rho^N_k}\rp-\xi\lp\frac{y'}{\rho^N_k}\rp\rv  
\leq|y-y'|.
\end{align*}
Since $\widehat\rho\Nk(0)=0$, we deduce
$ 
\lv\widehat \rho\Nk(x)\rv\leq \lv x\rv.
$\\
3. This upper estimate is a straightforward consequence of $\lv \xi\rv_{\infty}\leq 2$;
this concludes the proof.
\endproof
Let  $\lp X^{N,m}_.\rp_{1\leq m \leq M}$, $\lp \DeltaWm{.}\rp_{1\leq m \leq M}$ and  $\lp \DeltaBm{.}\rp_{1\leq m \leq M}$
be independent realizations of  $X^N$, $\DeltaW{}$ and $\DeltaB{}$ respectively.
In a similar way, we introduce the following random variables and random functions:
\begin{align}
\zeta^{N,m}_k:=& \rho^{N,m}_k:= \lv p\mk\rv\sqrt{C_0}\vee 1,\nonumber\\
\widehat\zeta_k^{N,m}(x):=&\zeta^{N,m}_k\xi\lp\frac{x}{\zeta^{N,m}_k}\rp, \quad
\widehat\rho_k^{N,m}(x):=\rho^{N,m}_k\xi\lp\frac{x}{\rho^{N,m}_k}\rp,\; x\in\R.\label{def-widehat-rho-m}
\end{align}
An argument similar to that used to prove Lemma \ref{lem-inv} yields the following:
\begin{lemma} \label{lem-widehat-rho-Nm}
The random functions
 $\widehat\rho\Nm\k(.)$ defined above satisfy the following properties:
\begin{enumerate}
\item $\widehat\rho\Nm\k$ is bounded by $2\lv\rho\Nm\k\rv$ and is 1-Lipschitz.
\item $\rho\Nm\k$ and  $\rho^N\k$ have the same distribution.
\end{enumerate}
\end{lemma}
We now describe the numerical scheme
\begin{definition}\textit{Initialization.}
At time $t=t_N$, set $Y^{N,i,I,M}_{t_N}:=\alpha^{i,I,M}_N.p_N:=P_N\Phi\lp X_{t_N}^N\rp$
and $\beta^{i,I,M}_N=0$ for all $i\in\lb1,\dots,I\rb$.\\
\textit{Induction}
Assume that an approximation $Y^{N,i,I,M}_{t_l} $ 
is built for $l=k+1,\dots,N$ and set 
$Y^{N,I,I,M,m}\tkk:=\widehat\rho^{N,m}_k\lp \alpha\IIM_{k+1}.p^m_{k+1} \rp$ its realization along the $m$th simulation. \\
We use backward induction in time and forward induction on $i$.
For $i=0$, let $\alpha^{0,I,M}_k=\beta^{0,I,M}_k=0$.
For $i=1,\dots,I$, the vector $\theta\iIM\k:=\lp\alpha^{i,I,M}_k,\sqrt h\beta^{i,I,M}_k\rp$ 
is defined by (forward) induction as the arg min in $(\alpha,\beta)$
of the quantity:
\begin{align}
\unM\summ&\lv Y^{N,I,I,M,m}\tkk -\alpha. \pmk
+h f\lp X^{N,m}\tk,\alpha^{i-1,I,M}_{k}.p^m_k ,\beta^{i-1,I,M}_{k}.p_k^m\rp\right.\nonumber\\
&\left.+\DeltaBmk g\lp X^{N,m}\tkk,Y^{N,I,I,M,m}\tkk\rp 
-\beta. \qmk \DeltaWm{k+1}\rv^2\label{eq-4}.
\end{align}
This minimization problem is similar to \eqref{eq-9}
replacing the expected value by an average over $M$ independent realizations.
Note that $\theta\iIM\k=\lp\alpha^{i,I,M}_k,\sqrt h\beta^{i,I,M}_k\rp$ is a random vector.
We finally set:
\begin{align}\label{4-4-bis}
Y\NIIM\tk:=\widehat \rho\Nk\lp \alpha\IIM\k .p\k\rp,
 \sqrt h Z\NIIM\tk:=\widehat \zeta\Nk\lp \sqrt h \beta\iIM\k.p\k\rp,
\end{align}
\end{definition}
The following theorem gives an upper estimate of the $L^2$ 
error beetween $\lp Y\NII_.,Z\NI_.\rp$ and $\lp Y\NIIM_.,Z\NIIM_.\rp$
in terms of $\lv \zeta^N_.\rv$ and $\lv \rho^N_.\rv$; it is the main
result of this section.
We recall that by \eqref{def-AMk} $
\AMk=\bigcap_{j=k}^{N-1} 
\left\lbrace \left\| V^M_j-Id\right\|\leq h ,
\|P^M_j-Id \|\leq h\right\rbrace\in\F^{W,M}_T\vee\F^{B,M}_{t\k,T}$.
For $k=1,\dots,N-1$ set
\begin{align}
\e\k:=&\E\|v\k v\k^*-Id\|_F^2\lp\E\lv\rho\Nk\rvd+\E\lv\zeta\Nk\rvd\rp
+\E\lc\lv v\k\rvd\lv p\kk\rvd\rc\E\lv\rho^N\kk\rvd
\nonumber\\
&+h^2\E\lc \lv v\k\rvd\lp1+\lv X^N\k\rvd+\lv p\k\rvd\E\lv\rho\Nk\rvd+\unh\lv p\k\rvd\E\lv \zeta\Nk\rvd\rp\rc\nonumber\\
&
+h\E\lc\lp\lv v\k\rvd+\lv w^p\k\rvd\rp
\lp1+\lv X^N\tkk\rvd+\lv p\kk\rvd \E\lv \rho^N\kk\rvd\rp \rc.
\label{def-ek}
\end{align}

Choosing $N$ and then $M$ large enough, the following
result gives the speed of convergence of the Monte
Carlo approximation scheme of $Y^{N,I,I}$ and 
$Z^{N,I}$.
\begin{theorem}\label{theorem-stepIII}
There exists a constant $C>0$ such that
for $h$ small enough, for any $k=0,\dots,N-1$
and $M\geq 1$:
\begin{align*}
\EM:=&\E\lv Y\NII\tk -Y\NIIM\tk\rv^2+h\sum_{j=k}^{N-1}\E\lv Z^{N,I}_{t_j}-Z\NIIM_{t_j}\rv^2\\
\leq& 16 \sum_{j=k}^{N-1}\E\lc\lp\lv\zeta\Nj\rv^2
+\lv\rho^N_j\rv^2\rp1_{\lc\AMk\rc^c}\rc
+Ch^{I-1}\sum_{j=k}^{N-1}
\lp h^2+h\E\lv\rho^N_{j+1}\rvd+\E\lv\rho\Nj\rvd+\E\lv\zeta\Nj\rvd\rp\\
&+\frac{C}{hM}\sum_{j=k}^{N-1}\e_j.
\end{align*}
\end{theorem}
\subsection{Proof of Theorem \ref{theorem-stepIII}}
Before we start the proof, let us recall some results on regression (i.e. orthogonal projections).
Let $v=(v^m)_{1\leq m\leq M}$ be a sequence of vectors in $\R^n$. 
Let use define the $n\times n$ matrix 
$V^M:=\unM \summ v^mv^{m*}$,
suppose that $V^M$ is invertible 
and denote by $\lambda_{\min}\lp V^M\rp$ its smallest eigenvalue.
\begin{lemma}\label{lem-regression}
Under the above hypotheses, we have the following results:
Let $(x^m,m=1,\dots,M)$ be a vector in $\R^M$.
\begin{enumerate}
\item There exists a unique $\R^n$ valued vector $\theta_x$ satisfying
$ \theta_x =\underset{\theta\in\R^n}{\arg\inf}|x-\theta. v|_M^2$
where $\theta.v$ denotes the vector
$\lp \sum_{i=1}^n\theta(i)v^m(i), m=1\,\dots, M\rp$.
\item Moreover, we have
$\theta_x=\unM \lp V^M\rp^{-1}\summ x^mv^m\in\R^n$
\item The map $x\mapsto \theta_x$ is linear from $\R^M$ to $\R^n$ and
$\lambda_{\min}(V^M)|\theta_x|^2\leq |\theta_x .v|_M^2\leq |x|_M^2$.
\end{enumerate}
\end{lemma}
The following lemma gives a first upper estimate of $\EM$.
\begin{lemma} For every $M$ and $k=0,\dots,N-1$,
we have the following upper estimate
\begin{align*}
\EM
\leq & 
\E\lc\lv \alpha^{I,I}_k- \alpha^{I,I,M}_k\rv^21_{\AMk}\rc
+h\sum_{j=k}^{N-1}\E\lc\lv \beta_j^{I}-\beta_j\IIM\rv^21_{\AMj}\rc\\
&+16\E\lc\lv\rho^N_k\rv^2 1_{\lc\AMk\rc^c}\rc
+16\sum_{j=k}^{N-1}\E\lc\lv\zeta\Nj\rv^2 1_{\lc\AMj\rc^c}\rc.
\end{align*}
\end{lemma}
This lemma should be compared with inequality (31) in \cite{Gobet}.
\proof
Using the decomposition of $Y\NII$, $Y\NIIM$, $Z\NI$ and $Z\NIIM$,
Lemma \ref{lem-inv} (1) ,
we deduce
\begin{align*}
\EM
=&
\E\lc\lv \widehat\rho^N_k\lp\alpha^{I,I}_k.p_k\rp-\widehat\rho^N_k\lp \alpha^{I,I,M}_k.p_k\rp\rv^2 \rc\\
&+h\sum_{j=k}^{N-1}\E\lc\lv \frac{1}{\sqrt h}\widehat\zeta\Nj\lp\sqrt h\beta_j^{I}.p_j\rp
-\frac{1}{\sqrt h}\widehat\zeta\Nj\lp\sqrt h\beta_j\IIM .p_j\rp\rv^2\rc.
\end{align*}
Using hte partition $A^M_k$, $\lp A^M_k\rp^c$ where
$A^M_k$ is defined by \eqref{def-AMk},
Cauchy-Schwarz's inequality, Lemma \ref{lem-inv}  
and the independence of $\lp\alpha\IIM\k,\beta\IIM_j,1_{\AMk}\rp$ and $p\k$
we deduce:
\begin{align*}
\EM\leq &
 \E\lc\lv \alpha^{I,I}_k.p_k- \alpha^{I,I,M}_k.p_k\rv^21_{\AMk}\rc
+h\sum_{j=k}^{N-1}\E\lc\lv \beta_j^{I}.p_j-\beta_j\IIM. p_j\rv^21_{\AMj}\rc\\
&+2\E\lc\lp\lv \widehat\rho^N_k\lp\alpha^{I,I}_k.p_k\rp\rv^2+\lv\widehat\rho^N_k\lp \alpha^{I,I,M}_k.p_k\rp\rv^2\rp1_{\lc\AMk\rc^c}\rc\\
&+2\sum_{j=k}^{N-1}\E\lc\lp\lv \widehat\zeta\Nj\lp\sqrt h\beta_j^{I}.p_j\rp\rv^2
+\lv\widehat\zeta\Nj\lp\sqrt h\beta_j\IIM. p_j\rp\rv^2\rp 1_{\lc\AMj\rc^c}\rc\\
\leq & 
\E\lc\lp \alpha^{I,I}_k-\alpha^{I,I,M}_k\rp^* p_kp\k^*\lp \alpha^{I,I}_k-\alpha^{I,I,M}_k\rp1_{\AMk}\rc\\
&+h\sum_{j=k}^{N-1}\E\lc\lp \beta_j^{I}-\beta_j\IIM\rp^* p_jp_j^*\lp \beta_j^{I}-\beta_j\IIM\rp1_{\AMj}\rc\\
&+2\E\lc 8\lv\rho^N_k\rv^21_{\lc\AMk\rc^c}\rc
+2\sum_{j=k}^{N-1}\E\lc 8\lv\zeta\Nj\rv^2 1_{\lc\AMj\rc^c}\rc\\
\leq & 
\E p_kp\k^*\E\lc\lv \alpha^{I,I}_k- \alpha^{I,I,M}_k\rv^21_{\AMk}\rc
+h\sum_{j=k}^{N-1}\E p_jp_j^*\E\lc\lv \beta_j^{I}-\beta_j\IIM\rv^21_{\AMj}\rc\\
&+16\E\lc\lv\rho^N_k\rv^21_{\lc\AMk\rc^c}\rc
+16\sum_{j=k}^{N-1}\E\lc\lv\zeta\Nj\rv^2 1_{\lc\AMj\rc^c}\rc.
\end{align*}
This concludes the proof.
\endproof
We now upper estimate $\lv \theta\IIM_k-\theta\II_k\rv^2$ on the event $\AMk$. This will be done in severals lemmas below.
By definition $\|V\Mk-I\|\leq h$  on $\AMk$ for any $k=1,\dots,N$.
Hence for $h\in(0,1)$
\begin{align}
1-h\leq \lambda_{\min}\lp V\Mk(\omega)\rp\quad \mbox{\rm on}\;\AMk.\label{lem-VMk}
\end{align}
\begin{lemma}\label{lem-linalg}For every $\alpha\in\R^n$ and $k=1,\dots,N$, we have
$
\unM\summ\lv \alpha .p\mk\rvd\leq \lv \alpha\rvd\left\| P\Mk\right\|.
$
\end{lemma}
\proof 
The definition of the Euclidian norm and of $P\Mk$ imply
\begin{align*}
\unM\summ\lv \alpha .p\mk\rvd
=\alpha^*\unM\summ p\mk\lp p\mk\rp^*\alpha
=\alpha^*P\Mk\alpha
\leq\left \|P\Mk\right\| \lv \alpha\rvd;
\end{align*}
this concludes the proof.
\endproof
For $i=0,\dots,I$, we introduce the vector 
$x\iIM\k:=\lp x\iImM\k\rp_{m=1,\dots,M}$ 
defined for $m=1,\dots,M$ by:
\begin{align}
x\iImM\k:=&\wrNm\kk\lp\alpha\IIM\kk .p^m\kk\rp
 +hf\lp X\Nm\k,\alpha\iIM\k .p\k^m,\beta\iIM\k .p\k^m\rp\nonumber\\
& +\DeltaBmk g\lp X\Nm\kk,\widehat\rho\Nm\kk\lp\alpha\IIM\kk .p^m\kk\rp\rp.
\label{eq-xiImMk}
\end{align}
Using Lemma \ref{lem-regression}, we can rewrite equation \eqref{eq-4} as
follows: 
\begin{align}
\theta\iIM\k=&\arg\inf_{\theta}\lv x\imIM\k-\theta. v\mk\rvd_M
= \unM\lp \VM_k \rp^{-1}\summ x\imImM\k v\mk.\label{eq-pour-corolaire}
\end{align}
We will need the following 
\begin{lemma}\label{lem-alpha-F}
For all $k=0,\dots,N-1$ and every $I$, the random variables
$\alpha\IIM_k$ are $\F^{W,M}_T\vee\F^{B,M}_{t\k,T} $ measurable.
\end{lemma}
\proof
The proof uses backward indution on $k$ and forward induction on $i$.\\
\textit{Initialization.} 
Let $k=N-1$.
By definition $\alpha^{0,I,M}_{N-1}=0$. Let $i\geq1$ and suppose $\alpha^{i-1,I,M}_{N-1}\in\F^{W,M}_T\vee\F^{B,M}_{t_{N-1},T} $. 
Using \eqref{def-vmk} (resp. \eqref{def-VMk}), we deduce that $v^m_{N-1}\in\F^{W,m}_{T}\vee\F^{B,m}_{t_{N-1},T}$
(resp. $V^M_{N-1}\in\F^{W,M}_{T}\vee\F^{B,M}_{t_{N-1},T}$).\\
Futhermore \eqref{eq-xiImMk} shows that $x\imImM_{N-1}\in\F^{W,M}_{T}\vee\F^{B,M}_{t_{N-1},T}$
and hence \eqref{eq-pour-corolaire} implies that $\alpha^{i,I,M}_{N-1}\in\F^{W,M}_T\vee\F^{B,M}_{t_{N-1},T} $.  \\
\textit{Induction.} 
Suppose that $\alpha\IIM\kk\in \F^{W,M}_T\vee\F^{B,M}_{t\kk,T}$; 
we will prove by forward induction on $i$ that $\alpha^{i,I,M}_{k}\in\F^{W,M}_T\vee\F^{B,M}_{t_{k},T}$ for $i=0,\dots,I$. \\
By definition $\alpha^{0,I,M}_{k}=0$. 
Suppose $\alpha^{i-1,I,M}_{k}\in\F^{W,M}_T\vee\F^{B,M}_{t_{k},T} $;
we prove that  $\alpha^{i,I,M}_{k}\in\F^{W,M}_T\vee\F^{B,M}_{t_{k},T} $ by similar arguments. Indeed, 
 \eqref{def-vmk} (resp. \eqref{def-VMk}) implies that $v^m_{k}\in\F^{W,m}_{T}\vee\F^{B,m}_{t_{k},T}$
(resp. $V^M_{k}\in\F^{W,M}_{T}\vee\F^{B,M}_{t_{k},T}$),
while \eqref{eq-xiImMk} (resp. \eqref{eq-pour-corolaire}) yields $x\imImM_{k}\in\F^{W,M}_{T}\vee\F^{B,M}_{t_{k},T}$
(resp. $\alpha^{i,I,M}_{k}\in\F^{W,M}_T\vee\F^{B,M}_{t_{k},T} $).
This concludes the proof.
\endproof
The following Lemma gives an inductive upper estimate of 
$\lv\theta\iiIM_k-\theta\iIM_k\rv^2$.
\begin{lemma}\label{lemme-pour-C-tilde}
There exists $\widetilde C>0$ such that for small $h$, for $k=0,\dots,N-1$ and  for $i=1,...,I-1$
$
\lv\theta\iiIM_k-\theta\iIM_k\rv^2\leq \widetilde Ch\lv\theta\iIM_k-\theta\imIM_k\rv^2 \quad on\; \AMk.
$
\end{lemma}
\proof Using \eqref{lem-VMk} and Lemma \ref{lem-regression} (4), we obtain on $\AMk$
\begin{align*}
  (1-h)\lv\theta\iiIM_k-\theta\iIM_k\rv^2
\leq  
\lambda_{\min}\lp V\Mk\rp\lv\theta\iiIM_k-\theta\iIM_k\rv^2
\leq
 \lv x\iIM_k-x^{i-1,I,M}_k\rv_M^2.
\end{align*}
Plugging equation \eqref{eq-xiImMk} and using the Lipschitz property \eqref{cf} of $f$,
 we deduce
\begin{align*}
  (1-h)\lv\theta\iiIM_k-\theta\iIM_k\rv^2
\leq  
\frac{h^2L_f}{M}\summ&\lp \lv\lp\alpha\iIM_k- \alpha\imIM_k\rp.p\mk\rv^2\right.\\
&\left.+\lv\lp\beta\iIM_k- \beta\imIM_k\rp.p\mk\rv^2\rp.
\end{align*}
 Lemma \ref{lem-linalg} and the inequality $\|P\Mk\|\leq 2$,
yield
\begin{align*}
  (1-h)\lv\theta\iiIM_k-\theta\iIM_k\rv^2
\leq&
 \lp \lv\alpha\iIM_k- \alpha\imIM_k\rv^2
+\lv\beta\iIM_k- \beta\imIM_k\rv^2\rp h^2L_f\left\|P\Mk\right\|\\
\leq & 2hL_f\lv\theta\iIM_k-\theta\imIM_k\rv^2.
\end{align*}
This concludes the proof.
\endproof
For $\theta=\lp \alpha,\sqrt h \beta\rp$ set
$F_k(\theta):=\arg\inf_{\theta^*}\lv x_k^{I,M}(\theta)-\theta^*.v\k\rvd$
where\\ 
$x\k\IM\lp\theta\rp:=\rho^{N,m}\kk\lp\alpha\IIM\kk.p^m\kk\rp
+hf\lp X^{N,m}\tk,\alpha.p^m\k,\beta.p^m\k\rp
+\DeltaBmk g\lp X^{N,m}\tkk,\hat\rho^{N,m}\kk\lp\alpha\IIM\kk.p^m\kk\rp\rp.$
\begin{lemma}\label{lem-Fk}
On $\AMk$, the map $F\k$ is Lipschitz
with a Lipschitz constant $2hL_f(1-h)^{-1}$.
\end{lemma}
\proof
Using \eqref{lem-VMk} and Lemma \ref{lem-regression} (3), we obtain on $\AMk$
\begin{align*}
(1-h)\lv F\k\lp \theta_1\rp-F\k\lp \theta_2\rp\rvd
\leq
\lambda_{\min}\lp V^M\k\rp\lv F\k\lp \theta_1\rp-F\k\lp \theta_2\rp\rvd
\leq
\lv x\k\IM\lp\theta_1\rp-x\k\IM\lp\theta_2\rp\rvd.
\end{align*}
Using the Lipschitz property \eqref{cf} of $f$,
Lemma \ref{lem-linalg} and the inequality $\|P\Mk\|\leq 2$,
we deduce that on $A^M_k$:
\begin{align*}
(1-h)\lv F\k\lp \theta_1\rp-F\k\lp \theta_2\rp\rvd
\leq&
\frac{h^2L_f}{M}\summ\lp \lv\alpha_1.p\mk- \alpha_2.p\mk\rv^2
+\lv\beta_1.p\mk- \beta_2.p\mk\rv^2\rp.\\
\leq&
 \lv\alpha_1- \alpha_1\rv^2h^2L_f\left\|P\Mk\right\|
+\lv\beta_1- \beta_2\rv^2h^2L_f\left\|P\Mk\right\|\\
\leq & 2hL_f\lv\theta_1-\theta_2\rv^2;
\end{align*}
this concludes the proof.
\endproof
The Lipschitz property of $F\k$ yields the following:
\begin{corollary}\label{cor-lem1}
(i) For $h$ small enough, on $\AMk$, there exists a unique random vector 
$\theta_k\inftyIM:=\lp \alpha_k\inftyIM,\sqrt h\beta\inftyIM\k\rp$ such that
\begin{align}
  \theta_k\inftyIM=\unM\lp \VM\k \rp^{-1}\summ x\inftyImM\k v\mk
=\arg\inf_{\theta}\lv x^{I,M}\k\lp\theta\inftyIM\k\rp
-\theta.v_k\rv^2_M,\label{eq-theta-infIM}
\end{align}
where for $\theta=\lp\alpha,\sqrt h\beta\rp$,
$x^{I,M}_k\lp\theta\rp:=\lp x^{I,m,M}_k\lp\theta\rp\rp_{m=1,\dots,M}$
denotes the vector
with components
\begin{align*}
x^{I,m,M}_k\lp\theta\rp:=&\wrNm\kk\lp\alpha\IIM\kk .p\mkk\rp
+hf\lp X\Nm\tk,\alpha.p\mk,\beta.p\mk\rp\\
&+\DeltaBmk g\lp X\Nm\tkk,\widehat\rho\Nm\kk\lp\alpha\IIM\kk. p\mkk\rp\rp.
\end{align*}
Let $x\inftyIM\k=\lp x^{\infty,I,m,M}_k\rp_{m=1,\dots,M}
=\lp x^{I,m,M}\lp\theta\inftyIM\k\rp\rp_{m=1,\dots,M}.$
\\
(ii) Moreover there exits a constant $C>0$ such that for small $h$ and
any $k=0,\dots,N-1$
$$
\lv \theta\inftyIM_k-\theta\IIM_k\rv^2\leq Ch^I\lv \theta\inftyIM_k\rv^2.
$$
\end{corollary}
\proof
(i) This is a consequence of Lemma \ref{lem-Fk} since 
$2hL_f(1-h)^{-1}<1$ for small $h$.\\
(ii)
An argument similar to that used to
prove Lemma \ref{lem-Fk} implies
that for $i=1,\dots,I$
\begin{align*}
  (1-h)\lv\theta\inftyIM_k-\theta\IIM_k\rv^2
\leq & 2hL_f\lv\theta\inftyIM_k-\theta\ImIM_k\rv^2
\end{align*}
Since $\theta^{0,I,M}_k=0$, 
we conclude the proof.
\endproof
The following result, similar to Lemma \ref{lem-1}, will be crucial in subsequent estimates.
It requires some additional argument compared with similar estimates in \cite{Gobet}.
\begin{lemma}\label{lem-avant-R}
Let $U^m\kk$ be a $\F^{W,M}_T\vee\F^{B,M}_{\tkk,T}$ measurable random variable.
Then we have
\begin{align*}
\E\lc1_{\AMk}
U^m\kk
\DeltaBmk\rc =0.
\end{align*}
\end{lemma}
\proof
Using \eqref{def-mathfrak-Aj} and \eqref{def-AMk} we deduce
\begin{align*}
\E\lp1_{\AMk}U^m\kk\DeltaBmk\rp
 =\E\lp1_{A^M\kk}U^m\kk\E\lp\DeltaBmk1_{\mathfrak A_k}\Big|\F^{W,M}_T\vee\F^{B,M}_{t\kk,T}\rp\rp
\end{align*}
Recall that 
$ \mathfrak A_k=\left\lbrace \| V^M_k-Id\|\leq h ,\|P^M_k-Id \|\leq h\rb$.
We will prove that
\begin{equation}\label{eq-f-symmetric}
1_{\mathfrak A_k}=f\lp\DeltaB{k}^1,\dots,\DeltaB{k}^M \rp
\end{equation}
with a symmetric function $f$, that is 
$f\lp \beta_1,\dots,\beta_M\rp=f\lp -\beta_1,\dots,-\beta_M\rp$
for any $\beta\in\R^M$.\\
Suppose at first that \eqref{eq-f-symmetric} is true.
Since the distribution of the vectors 
$\lp \DeltaB{k}^1,\dots,\DeltaB{k}^M \rp$ 
and $\lp-\DeltaB{k}^1,\dots,-\DeltaB{k}^M \rp$
are the same, the independence of 
$\lp \DeltaB{k}^l,l=1,\dots,M\rp$ 
and $\F^{W,M}_T\vee\F^{B,M}_{t\kk,T}$ yields
\begin{align*}
\E\lp\DeltaBmk1_{\mathfrak A_k}\Big|\F^{W,M}_T\vee\F^{B,M}\tkk\rp
=&
\E\lp\DeltaBmk f\lp\DeltaB{k}^1,\dots,\DeltaB{k}^M \rp\rp\\
=&
\E\lp -\DeltaBmk f\lp-\DeltaB{k}^1,\dots,-\DeltaB{k}^M \rp\rp\\
=&
-\E\lp \DeltaBmk f\lp\DeltaB{k}^1,\dots,\DeltaB{k}^M \rp\rp.
\end{align*}
Which concludes the proof.\\
Let us now prove \eqref{eq-f-symmetric}. 
Clearly, it is enough to prove to prove that each norm 
involved in the definition of $\mathfrak A\k$
is of this form.
Let $A$ be one of the matrices $ V^M_k$ or $P^M_k $. 
Now we will compute the characteristic polynomial 
$\chi$ of the matrix $A-Id$ and prove that its
coefficients are symmetric.

Let $q^m$ be $p^m_k$ or $v^m_k$.
We reorganize $q^m$ as 
$q^m=\lp q^m_1, q^m_2\DeltaBmk \rp ^*$,
where $q^m_1$ are the elements of $q^m$
independent of $\DeltaBmk$,
and $q^m_2$ is independent of $
\DeltaBmk$.
So we have 
$$
q^m\lp q^{m}\rp^* 
=
\lp
\begin{array}{c|c}
q^m_1\lp q^{m}_1\rp^* 
&
q^m_1\lp q^{m}_2\rp^* \DeltaBmk
\\
\hline
q^m_2\lp q^{m}_1\rp^*\DeltaBmk
&
q^m_2\lp q^{m}_2\rp^*\DeltaBmk
\end{array}
\rp
$$
Let
$
A =
\unM \summ
q^m\lp q^{m}\rp^*
$;
then the characteristic polynomial 
of the matrix
$A - Id$
is given by
$$
\chi\lp A-Id\rp(X)
=
\det \lp
\begin{array}{c|c}
B-(X+1)Id
&
C
\\
\hline
C^*
&
D-(X+1)Id
\end{array}
\rp
$$
where
\begin{align*}
B:= &
\unM\summ 
q^m_1q^{m,*}_1
\in M_{I_1\times I_1}\lp \R\rp,
\quad
C:= 
\unM\summ
q^m_1q^{m,*}_2\DeltaBmk
\in M_{I_1\times I_2}\lp\R\rp,
\\
D :=& \unM\summ q^m_2q^{m,*}_2\lv\DeltaBmk\rvd
\in M_{I_2\times I_2}\lp\R\rp.
\end{align*}
Set
$J_1=\lb 1,\dots,I_1\rb$ and
$J_2=\lb I_1+1,\dots,I_1+I_2\rb$,
and for $\sigma \in \mathfrak S_{I_1+I_2}$
the following sets
$
\mathcal H (\alpha,\sigma,\beta)
= 
\lb i\in J_{\alpha}, \sigma(i)\in J_{\beta} \rb,
$
for $\alpha,\beta\in\lbrace1,2\rbrace$.
Using the definition of the determinant, 
we have
\begin{align*}
\chi(A-Id)(X) 
=&  
\sum_{\sigma\in \mathfrak S_{I_1+I_2}}
\epsilon(\sigma)
\prod_{i\in \mathcal H(1,\sigma,1)}
\lc B(i,\sigma(i))-(X+1)\delta_{i, \sigma(i)} \rc
\\
&
\prod_{i \in \mathcal H(1,\sigma,1)}
C(i,\sigma(i))
\prod_{i \in \mathcal H(2,\sigma,1)}
C(\sigma(i),i)
\prod_{i \in \mathcal H(2,\sigma,2)}
\lc D(i,\sigma(i)) - (X + 1)\delta_{i, \sigma(i) } \rc
\end{align*}
Since we have the relation
$ \lv \mathcal H(1,\sigma, 1)\rv
+ \lv \mathcal H(1, \sigma, 2)\rv
=\lv J_1\rv= I_1
$
and 
$
\lv \mathcal H(1,\sigma, 1)\rv 
+
\lv \mathcal H(2,\sigma,1)\rv
=
\lv J_1\rv=I_1
$,
we deduce that
$\lv \mathcal H(1,\sigma,1)\rv
+
\lv \mathcal H(2,\sigma,1)\rv
$
is even.
Therefore, the power of $\DeltaBmk$
in $\chi(A-Id)(X)$ is even, which concludes the proof.
\endproof
As a corollary, we deduce the following identities
\begin{corollary}\label{lem-pour-R} For $k=0,\dots,N-1$, we have
\begin{align}
&\E\lc 1_{\AMk}\summ \wrNm\kk
\lp\alpha\IIM\kk .p\mkk\rp
\DeltaBmk g\lp X\Nm\tkk,\widehat\rho\Nm\kk\lp\alpha\IIM\kk. p\mkk\rp\rp\rc=0, \label{eq-4-21-1}\\
&\E\lc 1_{\AMk} \lp
 \wrNm\kk\lp\alpha\II\kk.p\mkk\rp
  -
 \wrNm\kk\lp\alpha\IIM\kk.p\mkk\rp
\rp\right .\nonumber\\
&\quad\quad\left .\DeltaBmk \lp g\lp X\Nm\tkk,\alpha\II\kk.p\mkk\rp 
- g\lp X\Nm\tkk,\widehat\rho\Nm\kk\lp\alpha\IIM\kk.p\mkk\rp\rp
\rp\rc =0\label{eq-4-21-2}
\end{align}
\end{corollary}
\proof
Indeed, $X^{N,m}\tkk\in\F^{W,M}_T$. 
Futhermore, \eqref{def-widehat-rho-m}, Lemma \ref{lem-alpha-F} 
and the definition of $p^m\kk$ imply that
$\widehat\rho\Nm\kk\lp\alpha\IIM\kk.p\mkk\rp, 
\widehat\rho\Nm\kk\lp\alpha\II\kk.p\mkk\rp 
\in \F^{W,M}_T\vee\F^{B,M}_{t\kk,T}$.
Thus Lemma \ref{lem-avant-R} concludes the proof.
\endproof
The following result provides an $L^2$ bound of 
$\theta\inftyIM\k$ in terms of $\rho^N\kk$.
\begin{lemma}\label{lem-bound-thetainfIM}
There exists a constant $C$ such that, for every $N$
and $k=0,\dots,N-1$,
$
\E\lc 1_{\AMk}\lv\theta\inftyIM\k \rv^2\rc
\leq 
C
\E\lv \rho^{N}_{k+1}\rv^2 +Ch.
$
\end{lemma}
\proof
Using \eqref{lem-VMk}, Lemma \ref{lem-regression} (3)
and Corollary \ref{cor-lem1} (i) we have on $A^M_k$
\begin{align*}
(1-h)\lv\theta\inftyIM\k \rv^2
\leq \lambda_{\min}(\VMk)\lv\theta\inftyIM\k\rv^2
\leq \lv x\inftyIM_k\rv_M^2.
\end{align*}
Using (N11), taking expectation, using Young's inequality 
and \eqref{eq-4-21-1},
we deduce for any $\e>0$,
$k=0,\dots,N-1$,
\begin{align*}
(1-h)\E\lc 1_{\AMk}\lv\theta\inftyIM\k \rv^2\rc
\leq &
\sum_{j=1}^3T^{I,M}\k(j),
\end{align*}
where
\begin{align*}
 T^{I,M}\k(1):=&
\frac{1}{M}\lp1+\frac{h}{\e}\rp
\summ\E \lc 1_{\AMk}\lv \widehat{\rho}^{N,m}_{k+1}
\lp\aIIMkk.\pmkk \rp\rv^2\rc,\\
T^{I,M}\k(2):=&
\frac{h^2}{M}\lp 1+2\frac{\e}{h}\rp 
\summ \E\lc 1_{\AMk}
\lv f\lp\XNm\tk,\aiIMk.\pmk,\biIMk.\qmk \rp \rv^2\rc,\\
T^{I,M}\k(3):=&
\frac{1}{M}\lp1+\frac{h}{\e}\rp
\summ\E\lc 1_{\AMk}\lv\DeltaBmk 
g\lp \XNm\tkk,\widehat\rho\Nm\kk\lp\alpha\IIM\kk.\pmkk\rp \rp \rv^2\rc.
\end{align*}
Lemma \ref{lem-widehat-rho-Nm} yields
\begin{align}
T^{I,M}\k(1)
\leq  4\frac{1}{M}\lp1+\frac{h}{\e}\rp\summ\E\lv \rho^{N,m}_{k+1}\rv^2
\leq  4\lp1+\frac{h}{\e}\rp\E\lv \rho^{N}_{k+1}\rv^2.\label{eqT1}
\end{align}
The Lipschitz condition \eqref{cf} of $f$,
Lemma \ref{lem-linalg} and the inequalities 
$\left\|P\Mk\right\|\leq2$ 
valid on $\AMk$ imply
\begin{align}
T^{I,M}\k(2)\leq &
 2L_fh(h+2\e)\unM\summ\E\lc 1_{\AMk}\lv\alpha\inftyIM\k.p\mk\rvd 
+ 1_{\AMk}\lv\beta\inftyIM\k.p\mk\rvd\rc\nonumber\\
&+2 h(h+2\e)\unM\summ\E\lv  f\lp\XNm\tk,0,0\rp\rv^2\nonumber\\
\leq & 
2L_fh(h+2\e)\E\left\lbrace 1_{\AMk}\lp\lv\alpha\inftyIM\k \rvd 
+\lv\beta\inftyIM\k \rvd\rp\|P^M\k\|\right\rbrace
+2 h(h+2\e)\E\lv  f\lp X^N\tk,0,0\rp\rv^2\nonumber\\
\leq & 
4L_fh(h+2\e)\E \lc1_{\AMk}\lp\lv\alpha\inftyIM\k \rvd
+\lv\beta\inftyIM\k \rvd\rp\rc
+2 h(h+2\e)\E\lv  f\lp X^N\tk,0,0\rp\rv^2.\label{eqT2}
\end{align}
Finally, since $\DeltaBmk$ is independent of 
$\F^{W,M}_T\vee\F^{B,M}_{t\kk,T}$
for every $m=1,\dots,M$,
the Lipschitz property \eqref{cg} of $g$ 
and Lemma \ref{lem-widehat-rho-Nm} (1) yield
for $m=1,\dots,M$
\begin{align*}
\E\lc 1_{\AMk}\right. &\left.\lv\DeltaBmk g\lp \XNm\tkk,\widehat\rho\Nm\kk\lp\alpha\IIM\kk.\pmkk\rp \rp \rv^2\rc\\
=&
\E\lc 1_{A^M\kk}1_{\mathfrak A\k}\lv\DeltaBmk g\lp \XNm\tkk,\widehat\rho\Nm\kk\lp\alpha\IIM\kk.\pmkk\rp \rp \rv^2\rc\\
=&
\E\lc 1_{A^M\kk}\lv g\lp \XNm\tkk,\widehat\rho\Nm\kk\lp\alpha\IIM\kk.\pmkk\rp \rp \rv^2
\E\lp 1_{\mathfrak A\k}\lv\DeltaBmk\rvd\big|\F^{W,M}_{t_N}\vee\F^{B,M}_{t\kk,T}\rp \rc\\
\leq&
h\E\lc 1_{A^M\kk}\lv 
g\lp \XNm\tkk,\widehat\rho\Nm\kk\lp\alpha\IIM\kk.\pmkk\rp \rp \rv^2\rc\\
&\leq 
8L_gh\E\lc 1_{A^M\kk}\lv \rho\Nm\kk  \rv^2\rc
+2h\E\lc 1_{A^M\kk}\lv g\lp \XNm\tkk,0\rp\rv^2\rc.
\end{align*}
Therefore,
\begin{align}
T^{I,M}\k(3)\leq &
8L_gh\lp1+\frac{h}{\e}\rp \E \lv \rho^N\kk  \rv^2 
+2h\lp1+\frac{h}{\e}\rp \E\lv g\lp X^N\tkk,0\rp\rv^2.\label{eqT3}
\end{align}
The inequalities \eqref{eqT1}-\eqref{eqT3} imply that for any $\e>0$
and $h\in(0,1]$,
\begin{align*}
(1-h)\E\lc 1_{\AMk}\lv\theta\inftyIM\k \rv^2\rc
\leq &
 \left\lbrace 4\lp1+\frac{h}{\e}\rp +8L_gh\lp1+\frac{h}{\e}\rp    \right\rbrace
\E\lv \rho^{N}_{k+1}\rv^2\\
&+4L_fh(h+2\e)\E\lc 1_{\AMk}\lp\lv\alpha\inftyIM\k \rvd
+\lv\beta\inftyIM\k \rvd\rp\rc\\
&+2 h(h+2\e)\E\lv  f\lp X^N\tk,0,0\rp\rv^2
+2h\lp1+\frac{h}{\e}\rp \E\lv g\lp X^N\tkk,0\rp\rv^2.
\end{align*}
Choose $\e$ such that $8L_f\e=\frac{1}{4}$ so that $4L_f(h+2\e)=\frac{1}{4}+4L_fh$.
For $h$ small enough (that is $h\leq \frac{1}{4(4L_f+\undeux)}$), we have
$4L_f(h+2\e)\leq \undeux(1-h)$. Hence,
we deduce
$
\undeux(1-h)\E 1_{\AMk}\lv\theta\inftyIM\k \rv^2
\leq 
C
\E\lv \rho^{N}_{k+1}\rv^2+Ch
$,
which concludes the proof.
\endproof
The next result yields an upper estimate of the 
$L^2$-norm of $\theta\IIM\k-\theta\II\k$ in terms 
of $\theta\inftyIM\k-\theta\inftyI\k$.
\begin{lemma}\label{lem-TIM-TI}
There is a constant $C$ such that 
for every $N$ large enough and all $k=0,\dots,N-1$,
$$
\E \lc 1_{\AMk} \lv\theta\IIM\k-\theta\II\k\rv^2\rc\leq(1+Ch)\E\lc 1_{\AMk}\lv \theta\inftyIM\k -\theta\inftyI\k\rv^2 \rc
+  Ch^{I-1} \lp \E\lv\rho\Nk\rvd+\E\lv\zeta\Nk\rvd\rp.
$$
\end{lemma}
\proof We decompose $\theta\IIM\k-\theta\II\k$ as follows:
$$
\theta\IIM\k-\theta\II\k=\lp \theta\inftyIM\k -\theta\inftyI\k\rp + \lp \theta\IIM\k-\theta\inftyIM\k\rp - \lp\theta\II\k - \theta\inftyI\k\rp.
$$
Young's inequality implies 
\begin{align*}
\lv\theta\IIM\k-\theta\II\k\rv^2
=&
(1+h)\lv \theta\inftyIM\k -\theta\inftyI\k\rv^2 
+ 2\lp 1+\unh\rp
\lp
\lv\theta\II\k - \theta\inftyI\k\rv^2
+ 
\lv \theta\IIM\k-\theta\inftyIM\k\rv^2 
\rp.
\end{align*}
Taking expectation over the set $\AMk$, using Lemma \ref{lem-33} and 
the fact that $\alpha^{i,I}\k$ 
and $\beta^I\k$ are deterministic, we deduce
\begin{align*}
\E\lc 1_{\AMk}\lv\theta\IIM\k-\theta\II\k\rv^2\rc
\leq&
(1+h)\E \lc 1_{\AMk}\lv \theta\inftyIM\k -\theta\inftyI\k\rv^2 \rc
+2\lp 1+\unh\rp L_f^Ih^{2I}\E\lv \rho\Nk\rvd
\\
& 
+ 2\lp 1+\unh\rp \E \lc 1_{\AMk}\lv \theta\IIM\k-\theta\inftyIM\k\rv^2\rc.
\end{align*}
Since $\theta\inftyI\k$ is deterministic,
Corollary \ref{cor-lem1} (ii) and again Lemma \ref{lem-33} yield
\begin{align*}
\E\lc 1_{\AMk}\lv\theta\IIM_k
- \theta\inftyIM_k\rv^2\rc
\leq &
Ch^I \E \lc 1_{\AMk} \lv \theta\inftyIM_k\rv^2\rc\\
\leq &
Ch^I \lp \E\lv\rho\Nk\rvd+\E\lv\zeta\Nk\rvd\rp+Ch^I \E \lc 1_{\AMk} \lv \theta\inftyIM_k-\theta\inftyI\k\rv^2\rc.
\end{align*}
Therefore, we deduce 
\begin{align*}
\E\lc 1_{\AMk}\lv\theta\IIM\k-\theta\II\k\rv^2\rc\leq &(1+h)\E \lc 1_{\AMk}\lv \theta\inftyIM\k -\theta\inftyI\k\rv^2\rc 
+ 2\lp 1+\unh\rp Ch^I    \lp \E\lv\rho\Nk\rvd+\E\lv\zeta\Nk\rvd\rp\\
&+ 2\lp 1+\unh\rp 
Ch^I \E \lc 1_{\AMk} \lv \theta\inftyIM_k-\theta\inftyI\k\rv^2\rc
 +2\lp 1+\unh\rp L_f^Ih^{2I}\E\lv \rho\Nk\rvd \\
\leq &(1+Ch)\E \lc 1_{\AMk}\lv \theta\inftyIM\k -\theta\inftyI\k\rv^2 
\rc
+  Ch^{I-1} \lp \E\lv\rho\Nk\rvd+\E\lv\zeta\Nk\rvd\rp,
\end{align*}
which concludes the proof.
\endproof
The rest of this section is devoted to upper estimate
$ \theta\inftyIM\k -\theta\inftyI\k$ on $\AMk$.
We at first decompose $\theta\inftyI\k-\theta\inftyIM\k$
as follows:
\begin{equation}\label{dec}
 \theta\inftyI\k -\theta\inftyIM\k=\sum_{i=1}^5\B_i,
\end{equation}
where $\B_2$, $\B_3$ and $\B_5$ introduce 
a Monte-Carlo approximation of some expected value
by an average over the $M$-realization:
for $k=0,\dots,N-1$,
\begin{align*}
\B_1:=&\lp Id-\VMkm\rp\theta\inftyI\k,\\
\B_2:=&\VMkm \lc\E \lp v_k\wrN\kk\lp \alpha\II\kk .p\kk\rp\rp 
-\frac{1}{M}\summ v\mk\wrNm\kk\lp \alpha\II\kk .p\mkk\rp\rc,\\
\B_3:=&
\VMkm h \lc\E\lp  v_k f\lp X\Nk,\alpha\inftyI\k .p\k,\beta^I\k .p\k\rp\rp 
- \unM\summ v\mk
f\lp X\Nm\tk,\alpha\inftyI_k.p\mk,\beta^I_k.p\mk\rp \rc,\\
\B_4:=&
\unM\lp \VM\k \rp^{-1}\summ v\mk \lc
 \wrNm\kk\lp\alpha\II\kk .p\mkk\rp
  -
 \wrNm\kk\lp\alpha\IIM\kk .p\mkk\rp
\right.
\\
&\qquad+ 
hf\lp X\Nm\tk,\alpha\inftyI_k.p\mk,\beta^I_k.p\mk\rp 
-hf\lp X\Nm\tk,\alpha\inftyIM_k.p\mk,\beta\inftyIM_k.p\mk\rp
 \\
&\qquad\left.+ 
\DeltaBmk \lc g\lp X\Nm\tkk,\widehat\rho\Nm\kk\lp\alpha\II\kk .p\mkk\rp\rp 
- g\lp X\Nm\tkk,\widehat\rho\Nm\kk\lp\alpha\IIM\kk .p\mkk\rp\rp\rc
\rc, \\
\B_5:=&
\VMkm \lc \E\lp v_k  \DeltaB{k} g\lp X^N\tkk,\alpha\II\kk .p\kk\rp\rp
-\unM\summ v\mk\DeltaBmk g\lp X\Nm\tkk,\widehat\rho\Nm\kk\lp\alpha\II\kk .p\mkk\rp\rp \rc.
\end{align*}
Note that compared to the similar decomposition in \cite{Gobet}, 
$\B_4$ is slightly different and  $\B_5$ is new.
Indeed, using equation \eqref{eq-34} and  \eqref{eq-theta-infIM}
and Lemma \ref{lem-inv} (1), we obtain:
\begin{align*}
\theta\inftyI\k-&\theta\inftyIM\k=\lp Id-\VMkm\rp\theta\inftyI\k+\VMkm\theta\inftyI\k-\theta\inftyIM\k\\
=&\B_1+\VMkm\E\lc v_k \lp\alpha\II\kk .p\kk+hf\lp X\Nk,\alpha\inftyI\k .p\k,\beta^I\k .p\k\rp 
+\DeltaB{k} g\lp X^N\tkk,\alpha\II\kk .p\kk\rp\rp\rc\\
&-\unM\lp \VM\k \rp^{-1}\summ v\mk\lc
 \wrNm\kk\lp\alpha\IIM\kk .p\mkk\rp
+hf\lp X\Nm\tk,\alpha\inftyIM_k.p\mk,\beta\inftyIM_k.p\mk\rp \right.\\
&\qquad\qquad\qquad\left.\qquad\qquad+\DeltaBmk g\lp X\Nm\tkk,\widehat\rho\Nm\kk\lp\alpha\IIM\kk .p\mkk\rp\rp \rc
\end{align*}

\begin{align*}
=& \sum_{ j\in\lbrace1,2,3,5\rbrace}\B_j
+\unM\lp \VM\k \rp^{-1}\lc\summ v\mk
 \wrNm\kk\lp\alpha\II\kk .p\mkk\rp  -
\summ v\mk
 \wrNm\kk\lp\alpha\IIM\kk .p\mkk\rp\rc
 \\
&+ \unM\lp \VM\k \rp^{-1}\summ v\mk\lc h
f\lp X\Nm\tk,\alpha\inftyI_k.p\mk,\beta^I_k.p\mk\rp 
-
hf\lp X\Nm\tk,\alpha\inftyIM_k.p\mk,\beta\inftyIM_k.p\mk\rp\rc\\
&+\unM\lp \VM\k \rp^{-1}\summ v\mk \DeltaBmk
\lc g\lp X\Nm\tkk,\widehat\rho\Nm\kk\lp\alpha\II\kk .p\mkk\rp \rp
- g\lp X\Nm\tkk,\widehat\rho\Nm\kk\lp\alpha\IIM\kk .p\mkk\rp\rp\rc
,
\end{align*}
which concludes the proof of \eqref{dec}.

The following lemmas provide upper bounds of the error terms  $\B_i$.
Recall that if $F$ is a matrix such that
 $\|Id-F\|<1$, then $F$ is inversible, $F^{-1}-Id=\sum_{k\geq1}(Id-F)^k$
 and 
\begin{align}
\|Id-F^{-1}\|\leq\frac{\|Id-F\|}{1-\|Id-F\|}\label{lem-classic}
\end{align}
Indeed,
$
F^{-1}=(Id-(Id-F))^{-1}=\sum_{k\geq0}(Id-F)^{k}
$
and
$
\|Id-F^{-1}\|\leq\sum_{k\geq1}\|(Id-F)^{k}\|
$.
\begin{lemma}
\label{lem-Xi}
(i)
Let $(U_1,...,U_M)$ be a sequence of iid centered random variables. 
Then we have
$
\E\lv\summ U_j\rv^2=M\E\lv U_1\rv^2
$.\\
(ii) We have
$\E \left\|\summ\lp v\mk \lp v\mk\rp^*
-Id\rp\right\|_F^2= M\E \| v\k  v\k^*-Id\|_F^2$.
\end{lemma}
\proof (i) The proof is straightforward.\\
(ii) Using (i) (N6) and (N7), we deduce
\begin{align*}
\E \left\|\summ\lc v\mk \lp v\mk\rp^*-Id\rc\right\|_F^2
=& \sum_{i,j}\E\left|\summ\lc v\mk \lp v\mk\rp^*-Id\rc(i,j)\right|^2\\
=& M\sum_{i,j}\E\left|\lc v\k \lp v\k\rp^*-Id\rc(i,j)\right|^2
= M\E \| v\k  v\k^*-Id\|_F^2;
\end{align*}
this concludes the proof of the Lemma.
\endproof
The following lemma provides a $L^2$ upper bound of $\B_1$.
Recall that $A^M_k$ is defined by \eqref{def-AMk}.
\begin{lemma}[Upper estimate of $\B_1$]\label{estim-B1}
There exist a constant $C$ such that for small $h$ and 
every $M\geq1$,
$$
\E\lc\lv \B_1\rv^21_{\AMk}\rc
\leq\frac{C}{M}\E\|v_kv_k^*-Id\|_F^2
\lp\E\lv \rho^N_k\rv^2+\E\lv\zeta^N\k\rvd\rp.
$$
\end{lemma}
\proof
On $\AMk$ we have $\|Id-V^M_k\|\leq h<1$;
and hence \eqref{lem-classic} implies 
$\|Id-(V^M_k)^{-1}\|
\leq\frac{\|Id-V^M_k\|}{1-\|Id-V^M_k\|} 
\leq\frac{\|Id-V^M_k\|}{1-h}$. 
Using the inequality $\|.\|\leq\|.\|_F$ we deduce
\begin{align*}
\E\lc\|Id-(V^M_k)^{-1}\|^21_{\AMk}\rc 
\leq \frac{1}{(1-h)^2}\E\lc 1_{\AMk}\left \|Id-V^M_k\right\|_F^2\rc.
\end{align*}
By definition $V^M_k=\unM\summ v\mk\lp v\mk\rp^*$; so using Lemma 
\ref{lem-Xi} we obtain
$
\E \lc 1_{\AMk}\left\|Id-V^M_k\right\|_F^2\rc
\leq 
\frac{1}{M}\E \| v\k  v\k^*-Id\|_F^2.
$
Therefore, since $\theta^{\infty,I}\k$ is deterministic, 
 Lemma \ref{lem-33} yields
\begin{align*}
  \E\lc\lv \B_1\rv^21_{\AMk}\rc 
\leq& \lv\theta^{\infty,I}\k\rvd\E\lc\| Id-(V^M_k)^{-1}\|^21_{\AMk}\rc\\
\leq &\frac{C}{M}\lp\E\lv \rho^N_k\rv^2+\E\lv\zeta^N\k\rvd\rp \E \|v\k  v\k^*-Id\|_F^2;
\end{align*}
this concludes the proof.
\endproof
The next lemma gives an upper bound of 
$\left\|\lp \VMk\rp^{-1}\right \|$ on $\AMk$.
\begin{lemma}\label{lem-bound-VMkm}
For  $h\in(0,\undeux)$, we have 
$
\|(V^M_k)^{-1}\|\leq 2 \quad on \; \AMk
$.
\end{lemma}
\proof
Using the triangular inequality and inequality 
\eqref{lem-classic}, we obtain on $\AMk$
\begin{align*}
\|(V^M_k)^{-1}\|\leq &\|Id\|+ \|Id-(V^M_k)^{-1}\|
\leq  1+ \frac{\|Id-V^M_k\|}{1-\|Id-V^M_k\|}
\leq  1+\frac{h}{1-h}
=  \frac{1}{1-h}.
\end{align*}
Since $h<\undeux$, the proof is complete.
\endproof
The following result provides an upper bound of $\B_2$.
This estimate should be compared with that given 
in \cite{Gobet} page 2192.
\begin{lemma}[Upper estimate of $\B_2$]\label{estim-B2}
There exists a constant $C>0$ such that for large $N$ 
and $k=0,\dots,N-1$,
$  \E\lc\lv \B_2\rv^21_{\AMk}\rc \leq  \frac{C}{M}\E\lc \lv  v_k\rvd\lv p\kk\rv^2\rc\E \lv\rho^N\kk\rv^2.$
\end{lemma}
\proof
We can rewrite $\B_2$ as follows:
\begin{align*}
\B_2
=&
-\frac{(V_k^M)^{-1}}{M}
\summ \lp v^m_k\widehat\rho^{N,m}_{k+1}\lp\alpha^{I,I}_{k+1}.p^m\kk\rp
-\E\lc v_k\widehat\rho^{N}_{k+1}\lp\alpha^{I,I}_{k+1}.p\kk\rp\rc\rp.
\end{align*}
Using Lemmas \ref{lem-bound-VMkm} and \ref{lem-Xi} (i), 
we obtain for small $h$
\begin{align*}
 \E&\lc\lv \B_2\rv^21_{\AMk}\rc
\leq  
\frac{4}{M^2}\E\lc 1_{\AMk}\lv\summ 
\lp v^m_k\widehat\rho^{N,m}_{k+1}\lp\alpha^{I,I}_{k+1}.p^m\kk\rp
-\E\lc v_k\widehat\rho^{N}_{k+1}\lp\alpha^{I,I}_{k+1}.p\kk\rp\rc\rp\rv^2\rc\\
\leq & 
\frac{4}{M}\E\lv  v_k\widehat\rho^{N}_{k+1}\lp\alpha^{I,I}_{k+1}.p\kk\rp
-\E\lc v_k\widehat\rho^{N}_{k+1}\lp\alpha^{I,I}_{k+1}.p\kk\rp\rc\rv^2
\leq  
\frac{4}{M}\E\lv  v_k\widehat\rho^{N}_{k+1}\lp\alpha^{I,I}_{k+1}.p\kk\rp
\rv^2.
\end{align*}
Using Lemma \ref{lem-inv} (2), Cauchy-Schwarz's inequality 
and Proposition \ref{a-priori-estimate}, since 
$\alpha^{I,I}\kk$ is deterministic
we deduce
\begin{align*}
 \E\lc\lv \B_2\rv^21_{\AMk}\rc
\leq &
 \frac{4}{M}\E\lc\lv  v_k\rvd\lv\alpha^{I,I}_{k+1}.p\kk
\rv^2\rc
\leq 
\frac{4}{M}\E\lc \lv  v_k\rvd\lv p\kk\rv^2\rc\E \lv\rho^N\kk\rv^2,
\end{align*}
which concludes the proof.
\endproof
The next lemma gives an upper estimate of the $L^2$-norm of $\B_3$
\begin{lemma}[Upper estimate of $\B_3$]\label{estim-B3}
There exists a constant $C$ such that 
for large $N$ and $k=0,\dots,N-1$,
\begin{align*}
\E\lc 1_{\AMk}\lv\B_3\rvd\rc
\leq&
C\frac{h^2}{M}\E\lc\lv v_k \rvd \lp 1+\lv  X\Nk\rvd+\lv p\k\rvd\E\lv \rho\Nk\rvd+\frac{1}{h}\lv p\k\rvd\E\lv \zeta\Nk\rvd\rp\rc
\end{align*}
\end{lemma}
\proof 
We take expectation on $\AMk$, use Lemmas \ref{lem-bound-VMkm} 
and \ref{lem-Xi} (i);
this yields for small $h$
\begin{align}
\E\lc 1_{\AMk}\lv\B_3\rvd\rc
\leq&
4\frac{h^2}{M}\E\lp\lv  v_k \rvd \lv f\lp X\Nk,\alpha\inftyI_k.p\k,\beta^I_k.p\k\rp \rvd\rp\label{EE1}
\end{align}
The Lipschitz condition \eqref{cf}, Cauchy-Schwarz's inequality and Proposition \ref{a-priori-estimate} imply
\begin{align*}
&\lv f\lp X\Nk,\alpha\inftyI_k.p\k,\beta^I_k.p\k\rp \rvd
\leq 
2L_f\lp \lv  X\Nk\rvd+\lv\alpha\inftyI_k.p\k\rvd+\lv\beta^I_k.p\k\rvd\rp
+2\lv f(0,0,0)\rvd\nonumber\\
\leq &
2L_f\lp \lv  X\Nk\rvd+\lv p\k\rvd\E\lv \rho\Nk\rvd+\frac{1}{h}\lv p\k\rvd\E\lv \zeta\Nk\rvd\rp
+2\lv f(0,0,0)\rvd,
\end{align*}
which together with \eqref{EE1} concludes the proof.
\endproof
The next result gives an upper estimate of $\B_4$ in $L^2$.
\begin{lemma}[Upper estimate of $\B_4$]\label{estim-B4}
Fix $\e>0$; there exist constants $C$ and $C(\e)$ such that
for $N$ large and $k=0,\dots,N-2$,
\begin{align*}
(1-h)\E\lc 1_{\AMk}\lv\B_4\rvd\rc
\leq & 
\lp 1+C(\e)h\rp 
\E\lc 1_{A^M\kk}\lv
 \alpha\II\kk -\alpha\IIM\kk 
\rvd\rc\\
&+ 
C\lp h+2\e\rp h\lp\E\lc 1_{\AMk}
\lv\alpha\inftyI_k-\alpha\inftyIM_k\rvd\rc
+\E\lc 1_{\AMk}\lv\beta^I_k -\beta\inftyIM_k\rvd\rc\rp.
\end{align*}
\end{lemma}
\proof
By definition, we have $\B_4=\unM\lp V\Mk\rp^{-1}\summ v\mk x_4^m$.
Let $x_4:=\lp x^m_4,m=1,\dots,M\rp$; then 
Lemma \ref{lem-regression} and inequality \eqref{lem-VMk} imply that on $\AMk$,
$(1-h)\lv\B_4\rvd\leq\lambda_{\min}(\VMk)\lv\B_4\rvd\leq \lv x_4\rvd_M$.
Taking expectation, using Young's inequality and \eqref{eq-4-21-2} in Corollary \ref{lem-pour-R},
we obtain for $\e>0$: 
$(1-h)\E\lc 1_{\AMk}\lv\B_4\rvd\rc
\leq \sum_{i=1}^3T_i $, where:
\begin{align*}
T_1
:= & 
\lp1+\frac{h}{\e}\rp\unM\summ\E\lc 1_{\AMk}\lv
 \wrNm\kk\lp\alpha\II\kk.p\mkk\rp
  -
 \wrNm\kk\lp\alpha\IIM\kk.p\mkk\rp
\rvd\rc,
\\
T_2:=& 
\lp 1+\frac{2\e}{h}\rp h^2\unM\summ\E\lc 1_{\AMk}\lv
f\lp X\Nm\tk,\alpha\inftyI_k.p\mk,\beta^I_k.p\mk\rp 
\right.\right.\\
&\qquad\qquad\qquad\qquad\qquad\qquad\left.\left.
-f\lp X\Nm\tk,\alpha\inftyIM_k.p\mk,\beta\inftyIM_k.p\mk\rp
\rvd\rc,
 \\
T_3:=&\lp1+\frac{h}{\e}\rp\unM\summ\E\lc 1_{\AMk}\lv 
\DeltaBmk\lc g\lp X\Nm\tkk,\widehat\rho\Nm\kk\lp\alpha\II\kk.p\mkk\rp \rp
\right.\right.\right.
\\
&\qquad\qquad\qquad\qquad\qquad\qquad
\left.\left.\left.
- g\lp X\Nm\tkk,\widehat\rho\Nm\kk\lp\alpha\IIM\kk.p\mkk\rp\rp\rc
\rvd\rc.
\end{align*}
Lemma \ref{lem-widehat-rho-Nm} (1) and Lemma \ref{lem-linalg} yield
\begin{align*}
T_1\leq &
\lp1+\frac{h}{\e}\rp\unM\summ\E \lc1_{\AMk}\lv
 \alpha\II\kk.p\mkk
  -
\alpha\IIM\kk.p\mkk
\rvd\rc\\
\leq &
\lp1+\frac{h}{\e}\rp\E\lc 1_{\AMk}\lv
\alpha\II\kk   -\alpha\IIM\kk \rvd \|P^M\kk\|\rc.
\end{align*}
Since $\AMk\subset A^M\kk$
and $\|P^M\kk\|\leq 1+h$ on $\AMk$, we deduce  
\begin{align}\label{majT1}
T_1\leq &
\lp1+\frac{h}{\e}\rp (1+h)\E\lc1_{A^M\kk}\lv
\alpha\II\kk   -\alpha\IIM\kk \rvd \rc.
\end{align}
Using property \eqref{cf}, Lemma \ref{lem-linalg} and
a similar argument, we obtain
for $0<h\leq1$:
\begin{align}
T_2&\leq
L_fh\lp h+2\e\rp \E\lc 1_{\AMk}\lp
\lv\alpha\inftyI_k-\alpha\inftyIM_k\rvd
+
\lv\beta^I_k-\beta\inftyIM_k\rvd
\rp\|P^M\k\|
\rc
\nonumber\\
&\leq
2L_fh\lp h+2\e\rp \E\lc 1_{\AMk}\lp
\lv\alpha\inftyI_k-\alpha\inftyIM_k\rvd
+\lv\beta^I_k-\beta\inftyIM_k\rvd
\rp\rc.\label{majT2}
\end{align}
Finally, since $\AMk=A^M\kk\cap\mathfrak A_k$ and $\DeltaBm{k}$
is independent of $\F^W\tk\vee\F^B_{t\kk,T}$,
we have using the Lipschitz property \eqref{cg}:
\begin{align*}
T_3&\leq
\lp1+\frac{h}{\e}\rp\unM\summ\E\lc 1_{A^M\kk}
\lv g\lp X\Nm\tkk,\widehat\rho\Nm\kk\lp\alpha\II\kk.p\mkk\rp \rp
- g\lp X\Nm\tkk,\widehat\rho\Nm\kk\lp\alpha\IIM\kk.p\mkk\rp\rp
\rvd \right.\\
&\qquad\qquad\qquad\qquad\qquad\qquad\left.\E\lp1_{\mathfrak A\k} \lv\DeltaBmk\rvd\big|\F^W_{t_N}\vee\F^B_{t\kk,T}\rp
\rc\\
&\leq
L_gh\lp1+\frac{h}{\e}\rp\unM\summ\E\lc 1_{A^M\kk}
\lv\widehat\rho\Nm\kk\lp\alpha\II\kk.p\mkk\rp
- \widehat\rho\Nm\kk\lp\alpha\IIM\kk.p\mkk\rp
\rvd\rc.
\end{align*}
So using Lemma again \ref{lem-widehat-rho-Nm} (1) 
and Lemma \ref{lem-linalg},
we deduce
\begin{align}
T_3
\leq&
L_gh\lp1+\frac{h}{\e}\rp(1+h)\E\lc 1_{A^M\kk}
\lv\alpha\II\kk 
- \alpha\IIM\kk 
\rvd\rc.\label{majT3}
\end{align}
The inequalities \eqref{majT1}-\eqref{majT3} 
conclude the proof.
\endproof
We decompose $v\k$ as $v\k=(v\k^o,v\k^p)$ where
$v\k^o$ contains all the elements in the basis
which are independent to $\DeltaB{k}$ and
$v^p\k=\DeltahB{k}w^p\k$.
with $w^p\k$ independent to $\DeltaB{k}$.
The next lemma gives an upper estimate of 
the $L^2$ norm of $\B_5$ on $\AMk$.
\begin{lemma}[Upper estimate of $\B_5$]\label{estim-B5}
There exists constant $C$ such that
for small $h$ and $k=0,\dots,N-1$,
\begin{align*}
\E\lc 1_{\AMk}\lv\B_5\rvd\rc
\leq&
\frac{Ch}{M}\E\lc\lp \lv v\k\rvd +\lv w^p\k\rvd\rp
\lp1+\lv X^N\kk\rvd+\lv p\kk\rvd \E\lv \rho^N\kk\rvd\rp \rc.
\end{align*}
\end{lemma}
\proof 
The proof is similar to that of Lemma \ref{estim-B3} which deals with $\B_3$.
Lemmas \ref{lem-bound-VMkm}, \ref{lem-Xi} and \ref{lem-inv} (1) yield
for small $h$
\begin{align*}
\E\lc 1_{\AMk}\lv\B_5\rvd\rc
\leq&
\frac{4}{M}\E\lv
 v\k\DeltaB{k} g\lp X^N\tkk,\widehat\rho^N\kk\lp\alpha\II\kk.p\kk\rp \rp
-\E\lc v\k\DeltaB{k} g\lp X^N\tkk,\widehat\rho^N\kk\lp\alpha\II\kk.p\kk\rp \rp\rc
\rvd\\
\leq&\frac{4}{M}\E\lv
 v\k\DeltaB{k} g\lp X^N\tkk,\alpha\II\kk.p\kk\rp 
\rvd.
\end{align*}
Then the decompostion of $v\k$ yields
\begin{align*}
\E\lc 1_{\AMk}\lv\B_5\rvd \rc
\leq&\frac{4}{M}\E\lv
 v\k^o\DeltaB{k} g\lp X^N\tkk,\alpha\II\kk.p\kk\rp \rvd
+\frac{4}{M}\E\lv
 v\k^p\DeltaB{k} g\lp X^N\tkk,\alpha\II\kk.p\kk\rp 
\rvd.
\end{align*}
Since $\DeltaB{k}$ is independent of $\F^W_T\vee\F^B_{t\kk,T}$, we deduce
\begin{align*}
\E\lc1_{\AMk}\lv\B_5\rvd \rc
\leq&\frac{4}{M}\E\lv
 v\k^o\DeltaB{k} g\lp X^N\tkk,\alpha\II\kk.p\kk\rp \rvd
+\frac{4}{M}\E\lv
 w\k^p\frac{\lv\DeltaB{k}\rv^2}{\sqrt h} g\lp X^N\tkk,\alpha\II\kk.p\kk\rp 
\rvd
\\
\leq&\frac{Ch}{M}\E\lv v\k^o g\lp X^N\tkk,\alpha\II\kk.p\kk\rp \rvd
+\frac{Ch}{M}\E\lv
 w\k^p g\lp X^N\tkk,\alpha\II\kk.p\kk\rp 
\rvd
\\
\leq&\frac{Ch}{M}\E\lc\lp\lv v\k\rvd+\lv w\k^p\rvd\rp
\lv g\lp X^N\tkk,\alpha\II\kk.p\kk\rp \rvd\rc.
\end{align*}
The Lipschitz condition \eqref{cg}, Cauchy-Schwarz's
and Young's inequalities together with  
Proposition \ref{a-priori-estimate} yield
\begin{align*}
  \lv g\lp X^N\tkk,\alpha\II\kk.p\kk\rp \rvd
\leq&2L_g\lp \lv X^N\tkk\rvd+\lv\alpha\II\kk\rvd\lv p\kk\rvd\rp+2 \lv g(0,0)\rvd\\
\leq&2L_g\lp \lv X^N\tkk\rvd+\lv p\kk\rvd\E\lv \rho^N\kk\rvd\rp+2 \lv g(0,0)\rvd.
\end{align*}
This concludes the proof.
\endproof
{\bf Final step of the proof of Theorem \ref{theorem-stepIII}}.
Young's inequality implies that for $h\in(0,1]$,
$
(b_1+b_2+b_3+b_4+b_5)^2
\leq \frac{8}{h} (b_1^2+b_2^2+b_3^2+b_5^2)+(1+h)b_4^2.
$
Recall that $\e\k$ has been defined in \eqref{def-ek}.
Then the decomposition \eqref{dec} and Lemmas 
\ref{estim-B1} and \ref{estim-B2}-\ref{estim-B5} 
yield for $\e>0$, small $h$ and $\e_k$ defined by \eqref{def-ek}:
\begin{align*}
\E&\lc 1_{\AMk}\lv \theta\inftyIM\k-\theta\inftyI\k\rvd\rc
\leq 
\frac{8}{h} \E\lc 1_{\AMk}
\sum_{j\in\lbrace 1,2,3,5\rbrace}\lv \B_j\rvd
\rc
+(1+h) \E\lc 1_{\AMk}\lv\B_4\rvd\rc\\
\leq &
\frac{C}{Mh}\e\k
+(1+C(\e)h)
\E\lc 1_{A^M\kk}\lv
 \alpha\II\kk -\alpha\IIM\kk 
\rvd\rc\\
&+(1+Ch) 
C\lp h+2\e\rp h\lp\E\lc 1_{\AMk}
\lv\alpha\inftyI_k-\alpha\inftyIM_k\rvd\rc
+\E\lc 1_{\AMk}\lv\beta^I_k -\beta\inftyIM_k\rvd\rc\rp,
\end{align*}
where in the last inequality, we have used Lemma \ref{lem-1-Ch}.
The definition of $\theta\inftyIM\k$ and $\theta\inftyI\k$, yield
for $h$ small enough:
\begin{align*}
&\lc1-(1+Ch)C(h+2\e)h\rc
\E\lc1_{\AMk}
\lv\alpha\inftyI_k-\alpha\inftyIM_k\rvd\rc
+h\E\lc 1_{\AMk}\lv\beta^I_k -\beta\inftyIM_k\rvd\rc\\
\leq &
\frac{C}{Mh}\e\k
+(1+C(\e)h)
\E\lc 1_{A^M\kk}\lv
 \alpha\II\kk -\alpha\IIM\kk 
\rvd\rc
+(1+Ch) 
C\lp h+2\e\rp h
\E\lc 1_{\AMk}\lv\beta^I_k -\beta\inftyIM_k\rvd\rc.
\end{align*}
Using again Lemma \ref{lem-1-Ch}, we obtain for some constant $C$ and $h$ small enough
\begin{align}
\E\lc 1_{\AMk}\lv \theta\inftyIM\k-\theta\inftyI\k\rvd\rc
\leq &
\frac{C}{Mh}\e\k
+(1+C(\e)h)
\E\lc 1_{A^M\kk}\lv
 \alpha\II\kk -\alpha\IIM\kk 
\rvd\rc\nonumber\\
&+(1+Ch) 
C\lp h+2\e\rp h
\E\lc 1_{\AMk}\lv\beta^I_k -\beta\inftyIM_k\rvd\rc.\label{ineq-5.31}
\end{align}
Using Corollary \ref{cor-lem1} (ii) 
and Lemma \ref{lem-bound-thetainfIM} we deduce
\begin{align}
\E\lc 1_{\AMk}\lv\beta^I_k -\beta\inftyIM_k\rvd\rc
\leq &
2\E\lc 1_{\AMk}\lv\beta^I_k -\beta\IIM_k\rvd\rc 
+\frac{2}{h}\E \lc 1_{\AMk}\lv\theta\IIM_k -\theta\inftyIM_k\rvd\rc\nonumber\\
\leq &
2\E\lc 1_{\AMk}\lv\beta^I_k -\beta\IIM_k\rvd\rc 
+Ch^{I-1}\E\lc 1_{\AMk}\lv\theta\inftyIM_k\rvd\rc\nonumber\\
\leq &
2\E\lc 1_{\AMk}\lv\beta^I_k -\beta\IIM_k\rvd\rc 
+Ch^{I-1}\E\lv \rho^{N}_{k+1}\rv^2 +Ch^I.\label{ineq-5.32}
\end{align}
Plugging \eqref{ineq-5.31} and \eqref{ineq-5.32} in Lemma \ref{lem-TIM-TI}, we obtain
for some constant $C$ and $h$ small enough
\begin{align*}
\E\lc1_{\AMk} \lv\theta\IIM\k-\theta\II\k\rv^2\rc
\leq &
\frac{C}{Mh}\e\k
+(1+C(\e)h)
\E\lc1_{A^M\kk}\lv
 \alpha\II\kk -\alpha\IIM\kk 
\rvd\rc \\
&
+ Ch^{I-1} \lp h^2+h\E\lv\rho^N\kk \rvd
+ \E\lv\rho\Nk\rvd+\E\lv\zeta\Nk\rvd\rp\\
&+(1+Ch) 
C\lp h+2\e\rp h
\E\lc1_{\AMk}\lv\beta^I_k -\beta\IIM_k\rvd\rc
\end{align*}
But 
$(1+Ch) 
C\lp h+2\e\rp=2\e C+h(C+C^2h+2\e C^2)
$
and we may choose $\e$ such that $2\e C=\undeux$,
so that
$
1-(1+Ch) C\lp h+2\e\rp=\undeux-(C+C^2h+\frac{C}{2})h
$. Using again Lemma \ref{lem-1-Ch} we obtain for
some constant $C$ and $h$ small enough:
\begin{align*}
&\E\lc1_{\AMk}\lv\alpha\II\k-\alpha\IIM\k\rvd\rc
+h\undeux\lp1-Ch\rp\E\lc1_{\AMk} \lv\beta\IIM\k-\beta^I\k\rv^2\rc\\
&\leq 
(1+Ch)\E\lc1_{\AMk}\lv\alpha\II\kk-\alpha\IIM\kk\rvd\rc
+C\frac{\e\k}{hM}
+ Ch^{I-1} \lp h^2+h\E\lv\rho^N\kk\rvd
+\E\lv\rho\Nk\rvd+\E\lv\zeta\Nk\rvd\rp
\end{align*}
So for small $h$,
\begin{align*}
&(1-Ch)\left\lbrace\E\lc1_{\AMk}\lv\alpha\II\k-\alpha\IIM\k\rvd\rc
+h\undeux\E\lc1_{\AMk} \lv\beta\IIM\k-\beta^I\k\rv^2\rc\right\rbrace\\
&\leq 
(1+Ch)\E\lc1_{\AMk}\lv\alpha\II\kk-\alpha\IIM\kk\rvd\rc
+C\frac{\e\k}{hM}
+ Ch^{I-1} \lp h^2+h\E\lv\rho^N\kk\rvd
+\E\lv\rho\Nk\rvd+\E\lv\zeta\Nk\rvd\rp
\end{align*}
Using the Lemma \ref{lem-1-Ch}, we obtain
\begin{align*}
&\E\lc1_{\AMk}\lv\alpha\II\k-\alpha\IIM\k\rvd\rc
+h\undeux\E\lc1_{\AMk} \lv\beta\IIM\k-\beta^I\k\rv^2\rc\\
&\leq 
(1+Ch)\E\lc1_{\AMk}\lv\alpha\II\kk-\alpha\IIM\kk\rvd\rc
+C\frac{\e\k}{hM}
+ Ch^{I-1} \lp h^2+h\E\lv\rho^N\kk\rvd
+\E\lv\rho\Nk\rvd+\E\lv\zeta\Nk\rvd\rp
\end{align*}

The Gronwall Lemma \ref{lem-Gronwall}
applied with $a_k=\E\lc1_{\AMk}\lv\alpha\II\k-\alpha\IIM\k\rvd\rc$\\
and $c_k=h\undeux\E\lc1_{\AMk} \lv\beta\IIM\k-\beta^I\k\rv^2\rc$
and the fact that $\alpha^{I,I,M}_N=\alpha^{I,I}_N$ concludes the proof.
\subsection*{Acknowledgments:} 
The author wishes to thank Annie Millet for
many helpful comments.


\begin{thebibliography}{99}
\bibitem{Aboura} 
Aboura O.,
On the discretization of backward doubly stochastic differential equations
\textit{Arxiv 0907.1406}

\bibitem{Aman} 
Aman A.,
Numerical scheme for backward doubly stochastic differential
equations
\textit{Arxiv 0907.2035}

\bibitem{Aman2} 
Aman A.,
Numerical scheme for backward doubly stochastic differential
equations
\textit{Arxiv 1011.6170}

\bibitem{BallyDiscretization} 
Bally V., Approximation scheme for solutions of BSDE
\textit{ Backward stochastic differential equations (Paris, 1995--1996),
  177--191, Pitman Res. Notes Math. Ser.}, 364, Longman, Harlow, 1997.

\bibitem{BouchardTouzi} Bouchard B., Touzi N.,
Discrete time approximation 
and Monte Carlo simulation of backward stochastic differential equations
\textit{Stochastic process and applications} \textbf{111} (2004) 175-206

\bibitem{ChevanceThesis} 
Chevance D., 
Numerical methods for backward stochastic differential equations
\textit{Numerical methods in finance,  232--244, 
Publ. Newton Inst., Cambridge Univ. Press, Cambridge,} 1997. 

\bibitem{Gobet} 
Gobet E., Lemor J.P., Warin X., 
A regression-based Monte Carlo method 
to solve Backward stochastic differential equations
\textit{The Annals of Applied Probability} 
2005, Vol. 15, No. 3

\bibitem{Kloeden} 
Kloeden P. E., Platen E.,
\textit{Numerical Solution of Stochastic Differential Equations}
Springer

\bibitem{MaYong}
Ma J., Yong J., 
\textit{Forward-Backward Stochastic Differential Equations 
and their Applications} Lecture Note in Math. 1702 Springer 1999

\bibitem{NualartPardoux anticipating}
Nualart D., Pardoux E.,
Stochastic calculus with anticipatng integrands 
\textit{Probability theory related fields} 78, 535-581 (1988)

\bibitem{PPquasilinear} 
Pardoux E., Peng S., Backward stochastic differential
equation and quasilinear parabolic partial differential equations.
In: {\sl B. L.Rozovski, R. B. Sowers (eds). Stochastic partial
equations and their applications.  Lect. Notes control Inf. Sci.}
{\bf 176}, $200-217$, Springer, Berlin, $(1992)$.

\bibitem{PPspde} 
Pardoux E., Peng S.,
Backward doubly stochastic differential equations 
and systems of quasilinear SPDEs
\textit{Probability Theory and Related Fields} (1994) 209-227
\bibitem{Zhang} Zhang J.,
A numerical scheme for BSDEs
\textit{The Annals of Applied Probability} Vol. 14 No. 1 (2004) 459-488
\end{thebibliography}
\end{document}